\documentclass[10pt,reqno,oneside]{amsart}
\usepackage{amsmath,amsthm,amssymb,amscd,amsfonts,mathrsfs,amsxtra,upref}
\usepackage[mathscr]{eucal}
\numberwithin{equation}{section}
\makeatletter

\def\th@plain{%
\let\thm@indent\noindent
\thm@headfont{\caps}
\let\thmhead\thmhead@plain
\let\swappedhead\swappedhead@plain
\thm@preskip.5\baselineskip\@plus.2\baselineskip\@minus.2\baselineskip
\thm@postskip\thm@preskip
\slshape
}

\def\th@remark{%
\let\thm@indent\noindent
\thm@headfont{\caps}
\let\thmhead\thmhead@plain
\let\swappedhead\swappedhead@plain
\thm@preskip.5\baselineskip\@plus.2\baselineskip\@minus.2\baselineskip
\thm@postskip\thm@preskip
\upshape
}

\def\th@boldremark{%
\let\thm@indent\noindent
\thm@headfont{\bfseries}
\let\thmhead\thmhead@plain
\let\swappedhead\swappedhead@plain
\thm@preskip.5\baselineskip\@plus.2\baselineskip\@minus.2\baselineskip
\thm@postskip\thm@preskip
\upshape
}

\makeatother

\theoremstyle{plain}
\newtheorem{Theorem}{Theorem}[section]
\newtheorem{Theorem--Definition}[Theorem]{Theorem--Definition}
\newtheorem{Corollary}[Theorem]{Corollary}
\newtheorem{Lemma}[Theorem]{Lemma}
\newtheorem{Lemma--Definition}[Theorem]{Lemma--Definition}

\newtheorem{Proposition--Definition}[Theorem]{Proposition--Definition}
\theoremstyle{remark}
\newtheorem{Remark}[Theorem]{Remark}%
\newtheorem{Remark--Definition}[Theorem]{Remark--Definition}%
\newtheorem{Definition}[Theorem]{Definition}

\newtheorem{Notation}[Theorem]{Notation}

\DeclareMathOperator{\at}{\text{@}}

\DeclareMathOperator{\Aut}{{Aut}}
\DeclareMathOperator{\card}{\#}

\DeclareMathOperator*{\const}{const}

\DeclareMathOperator{\dcr}{DCR}

\DeclareMathOperator{\Gcd}{G.C.D.}

\DeclareMathOperator{\Id}{Id}

\DeclareMathOperator{\supp}{supp}

\font\caps=cmcsc10

\font\twelvegtc=eufm10 scaled 1200

\font\ninegtc=eufm9
\font\sevengtc=eufm7

\newfam\gtcfam

\textfont\gtcfam=\twelvegtc
\scriptfont\gtcfam=\ninegtc
\scriptscriptfont\gtcfam=\sevengtc

\def\Bbb{\mathbb}

\def\Cal{\mathcal}

\def\Def{\overset{\text{def}}{=\!=}}

\begin{document}
\author[Y. Feler]{Yoel Feler}
\address{Technion\\ Haifa\\   Israel 32000}
\email{feler$\at$techunix.technion.ac.il}
\thanks{This material is based upon a PhD thesis supported by Technion}

\title{Spaces of geometrically generic configurations}

\begin{abstract}
\noindent Let $X$ denote either $\Bbb{CP}^m$ or $\Bbb{C}^m$.
We study certain analytic properties of the space $\Cal{E}^n(X,gp)$
of ordered geometrically generic $n$-point
configurations in $X$.
This space consists of all $q=(q_1,...,q_n)\in X^n$ such that
no $m+1$ of the points $q_1,...,q_n$ belong to a hyperplane in $X$.
In particular, we show that for $n$ big enough any holomorphic
map $f\colon\Cal{E}^n(\Bbb{CP}^m,gp)\to\Cal{E}^n(\Bbb{CP}^m,gp)$
commuting with the natural action of the symmetric group
$\mathbf{S}(n)$ in $\Cal{E}^n(\Bbb{CP}^m,gp)$
is of the form $f(q)=\tau(q)q=(\tau(q)q_1,...,\tau(q)q_n)$,
$q\in \Cal{E}^n(\Bbb{CP}^m,gp)$,
where $\tau\colon\Cal{E}^n(\Bbb{CP}^m,gp)
\to{\mathbf{PSL}(m+1,\mathbb{C})}$ is an $\mathbf{S}(n)$-invariant
holomorphic map. A similar result holds true for mappings of
the configuration space $\Cal{E}^n(\Bbb{C}^m,gp)$.
\end{abstract}

\maketitle
\tableofcontents

\section{Introduction}\label{Sec: Introduction}

\noindent In this paper we study certain analytic properties of
the spaces of geometrically generic point configurations
in projective and affine spaces.
\vskip5pt

\noindent The most traditional non-ordered
configuration space ${\mathcal C}^n={\mathcal C}^n(X)$
of a complex space $X$ consists of all $n$ point subsets
({\em ``configurations"}) $Q=\{q_1,...,q_n\}\subset X$.
If $X$ carries an additional geometric structure,
it may be taken into account. Say if $X$ is either the projective
space ${\mathbb{CP}}^m$ or the affine space ${\mathbb C}^m$
and $n>m$ then the space ${\mathcal C}^n={\mathcal C}^n(X,gp)$
of {\em geometrically generic configurations} consists of all
$n$ point subsets $Q\subset X$ such that no hyperplane in $X$
contains more than $m$ points of $Q$. The corresponding
{\em ordered} configuration space 
${\mathcal E}^n={\mathcal E}^n(X,gp)$ consists of all
$q=(q_1,...,q_n)\in X^n$ such that the set
$Q=\{q_1,...,q_n\}\subset X$ belongs to ${\mathcal C}^n(X,gp)$. 
The left action of the symmetric group ${\mathbf S}(n)$
on ${\mathcal E}^n$ defined by $\sigma (q_1,...,q_n)
= (q_{\sigma^{-1}(1)},...,q_{\sigma^{-1}(n)})$ is free;
the corresponding orbit map can be identified with
the tautological covering map
$p\colon{\mathcal E}^n\ni (q_1,...,q_n)
\mapsto\{q_1,...,q_n\}\in{\mathcal C}^n$, which is an
unbranched Galois covering with the Galois group
${\mathbf S}(n)$.
\vskip5pt

\noindent When $m=1$, that is, $X={\mathbb{CP}^1}$ or $X={\mathbb C}^1$,
the spaces of geometrically generic configurations in $X$ coincide
with the usual configuration spaces ${\mathcal C}^n(X)$
and ${\mathcal E}^n(X)$. But for $m>1$ the spaces ${\mathcal C}^n(X,gp)$
and ${\mathcal E}^n(X,gp)$ seem even more natural since they both are
Stein manifolds, unlike ${\mathcal C}^n(X)$ and ${\mathcal E}^n(X)$.
In the same time, the spaces ${\mathcal C}^n(X,gp)$ and ${\mathcal E}^n(X,gp)$
have more complicated topology; for instance, their fundamental groups
are highly nontrivial, whereas $\pi_1({\mathcal C}^n(X))
={\mathbf S}(n)$ and ${\mathcal E}^n(X)$ is simply connected.
\vskip5pt

\noindent Let $\Aut_g X$ be the subgroup of the holomorphic
automorphism group $\Aut X$ of $X$ consisting of all
elements $A\in\Aut X$ that respect the geometrical structure
of $X$. I. e., for $X={\mathbb{CP}}^m$ we have
$\Aut_g X=\mathbf{PSL}(m+1,\mathbb{C})=\Aut X$,
whereas for $X={\mathbb C}^m$ the group
$\Aut_g X=\mathbf{Aff}(m,\mathbb{C})$ of all
affine transformations of $\mathbb{C}^m$ is much smaller than
the huge group $\Aut(\mathbb{C}^m)$ of all biholomolomorphic
(or polynomial) automorphisms of $\mathbb{C}^m$. Notice that
in both cases $\Aut_g X$ is a complex Lie group.
The natural action of $\Aut_g X$ in $X$ induces
the diagonal $\Aut_g X$ actions in ${\mathcal C}^n$
and ${\mathcal E}^n$ defined by
\begin{equation*}
AQ=A\{q_1,...,q_n\}\Def\{A q_1,...,A q_n\}\ \ \forall\,
Q=\{q_1,...,q_n\}\in {\mathcal C}^n
\end{equation*}
and
\begin{equation*}
Aq=A(q_1,...,q_n)\Def(A q_1,...,A q_n)\ \ \forall\,
q=(q_1,...,q_n)\in {\mathcal E}^n\,.
\end{equation*}
Therefore, any automorphism $T\in\Aut_g X$
gives rise to the holomorphic endomorphism (i.e., holomorphic
self-mapping) $F_T$ of ${\mathcal C}^n$ defined by
$F_T(Q)=TQ=\{Tq_1,...,Tq_n\}$. This example may be generalized
by making $T$ depending analytically on a configuration
$Z\in{\mathcal C}^n$; a rigorous definition looks as follows.

\begin{Definition}\label{Def: tame map}
A holomorphic endomorphism $F$ of ${\mathcal C}^n$
is said to be {\em tame} if there is
a holomorphic map $T\colon {\mathcal C}^n\to \Aut_g X$ such that
$F(Q)=F_T(Q)\Def T(Q)Q$ for all $Q\in {\mathcal C}^n$.
Similarly, a holomorphic endomorphism $f$ of ${\mathcal E}^n$
is called {\em tame} if there are an $\mathbf{S}(n)$-invariant
holomorphic map $\tau\colon {\mathcal E}^n\to \Aut_g X$
and a permutation $\sigma \in \mathbf{S}(n)$
such that $f(q)=f_{\tau,\sigma}(q)\Def\sigma\tau(q)q$ for all
$q\in {\mathcal E}^n$.

When $X={\mathbb C}^m$, a holomorphic endomorphism $F$ of
${\mathcal C}^n({\mathbb C}^m, gp)$
is said to be {\em quasitame} if there is a holomorphic map
$T\colon {\mathcal C}^n({\mathbb C}^m, gp)
\to\mathbf{PSL}(m+1,\mathbb{C})$ such that
$F(Q)=F_T(Q)\Def T(Q)Q$ for all
$Q\in{\mathcal C}^n({\mathbb C}^m,gp)$.\footnote{Notice that
the latter condition implies that $T(Q)Q\subset{\mathbb C}^m$
for any $Q\in{\mathcal C}^n({\mathbb C}^m,gp)$.}
Similarly, a holomorphic endomorphism
$f$ of ${\mathcal E}^n={\mathcal E}^n({\mathbb C}^m, gp)$
is {\em quasitame} if there are an $\mathbf{S}(n)$-invariant
holomorphic map $\tau\colon{\mathcal E}^n({\mathbb C}^m,gp)
\to\mathbf{PSL}(m+1,\mathbb{C})$ and a permutation
$\sigma\in\mathbf{S}(n)$ such that
$\tau(q)q_1,...,\tau(q)q_n\in{\mathbb C}^m$
and $f(q)=\sigma\tau(q)q$ for any
$q=(q_1,...,q_n)\in {\mathcal E}^n({\mathbb C}^m,gp)$.
\hfill $\bigcirc$
\end{Definition}

\noindent The left ${\mathbf S}(n)$ action
on $\Cal{E}^n(X,gp)$ induces the left ${\mathbf S}(n)$
action on the set of all mappings
$f\colon\Cal{E}^n(X,gp)\rightarrow \Cal{E}^n(X,gp)$ defined by
$$
\sigma f=\sigma(f_1,...,f_n) =
(f_{\sigma^{-1}(1)},...,f_{\sigma^{-1}(n)})
$$
for $f=(f_1,...,f_n)\colon\Cal{E}^n(X,gp)
\rightarrow \Cal{E}^n(X,gp)$ and $\sigma\in{\mathbf S}(n)$.

\begin{Definition}\label{Def: strict equivariance}
A continuous map $f\colon{\mathcal E}^n(X,gp)
\to{\mathcal E}^n(X,gp)$
is said to be {\em strictly equivariant} if there exists
an automorphism $\alpha$ of the group ${\mathbf S}(n)$
such that
\begin{equation*}
f(\sigma q)=\alpha(\sigma)f(q) \ \
\text{for all} \ q\in{\mathcal E}^n(X) \
\text{and} \ \sigma\in{\mathbf S}(n)\,.
\end{equation*}
\hfill $\bigcirc$
\end{Definition}

\begin{Remark}
\label{Rmk: tame and quasitame are strictly equivariant}
Any tame endomorphism $f=f_{\tau,\sigma}$ of
${\mathcal E}^n$ is strictly equivariant;
the corresponding automorphism $\alpha\in\Aut{\mathbf S}(n)$
is just the inner automorphism $s\mapsto\sigma s\sigma^{-1}$.
The same holds true for quasitame maps. \hfill $\bigcirc$
\end{Remark}

\noindent The following theorem contains the
main results of this work.

\begin{Theorem}\label{Thm: endomorphisms of En(X)}
Let $m>1$, $n\ge m+3$  and $n\ne 2m+2$.

$a)$ Any strictly equivariant holomorphic map
$f\colon\Cal{E}^n(\mathbb{CP}^m,gp)\rightarrow
\Cal{E}^n(\mathbb{CP}^m,gp)$ is tame.

$b)$ Any strictly equivariant holomorphic map
$f\colon\Cal{E}^n(\mathbb{C}^m,gp)\rightarrow
\Cal{E}^n(\mathbb{C}^m,gp)$ is quasitame.
\vskip5pt

In more detail, for $m>1$,  $n\ge m+3$, $n\ne 2m+2$ and
a strictly equivariant holomorphic map
$f\colon\Cal{E}^n(X,gp)\rightarrow\Cal{E}^n(X,gp)$,
there exist an $\mathbf{S}(n)$-invariant holomorphic map
$\tau\colon\Cal{E}^n(X,gp) \to \mathbf{PSL}(m+1,\mathbb{C})$
and a permutation $\sigma\in\mathbf{S}(n)$ such that
\begin{equation}\label{eq: tame representation of f}
f(q)=\sigma\tau(q)q= (\tau(q)q_{\sigma^{-1}(1)},...,
\tau(q)q_{\sigma^{-1}(n)})
\end{equation}
for all $q=\left(q_1,\dots,q_n\right)\in\Cal{E}^n(X,gp)$.
\end{Theorem}

\begin{Remark}
\label{rm: map commuting with action of symmetric group} For a
holomorphic endomorphism $f$ of $\Cal{E}^n(X,gp)$ commuting with
the $\mathbf{S}(n)$-action the representation \eqref{eq: tame
representation of f} can be simplified as follows: there is an
$\mathbf{S}(n)$-invariant holomorphic map
$\tau\colon\Cal{E}^n(X,gp) \to \mathbf{PSL}(m+1,\mathbb{C})$ such
that
$$
f(q)=\tau(q)q \ \ \text{for all} \ q\in\Cal{E}^n(X,gp)\,.
$$
\hfill $\bigcirc$
\end{Remark}

\begin{Remark}\label{Rmk: aaaa}
The holomorphic endomorphisms of the spaces
${\mathcal C}^n({\mathbb{CP}}^1)$ and ${\mathcal C}^n({\mathbb C})$
were completely described by V. Lin
(see, for instance, \cite{Lin72b}, \cite{Lin79} and
\cite{LinSphere}). Moreover, due to the works of V. Zinde
\cite{Zinde74}-\cite{Zinde78} and the author \cite{Feler04a},
such a description is known now for holomorphic endomorphisms
of traditional configuration spaces of all non-hyperbolic algebraic
curves $\Gamma$. In all quoted papers, the first important step
includes a purely algebraic investigation of the fundamental groups
$\pi_1({\mathcal C}^n(\Gamma))$ and $\pi_1({\mathcal E}^n(\Gamma))$,
which are the braid group and the pure braid group of
$\Gamma$, respectively. Namely, it is shown that for $n>4$
any endomorphism of $\pi_1({\mathcal C}^n(\Gamma))$ with
a non-abelian image preserves the subgroup
$\pi_1({\mathcal E}^n(\Gamma))$. This property ensures
the lifting of any ``sufficiently non-trivial"
holomorphic self-map $F$ of ${\mathcal C}^n(\Gamma)$
to a strictly equivariant holomorphic self-map $f$
of ${\mathcal E}^n(\Gamma)$, which, in turn, can
be studied via certain analytic and combinatorial methods.

My original target in the present work was to use a similar
machinery in order to study the holomorphic endomorphisms
of the spaces ${\mathcal C}^n({\mathbb{CP}}^m,gp)$ and
${\mathcal C}^n({\mathbb C}^m,gp)$ for $m>1$.
I did not succeed however in investigating of the
corresponding fundamental groups and therefore restricted myself
to the study of strictly equivariant holomorphic endomorphisms
of the corresponding ordered configuration spaces.
\hfill $\bigcirc$
\end{Remark}

\begin{Remark}
Theorem \ref{Thm: endomorphisms of En(X)}$(b)$
is not complete, since at the moment I do not know whether
there are strictly equivariant holomorphic endomorphisms of
${\mathcal E}^n({\mathbb C}^m,gp)$ that are quasitame
but not tame.

Although I think that Theorem \ref{Thm: endomorphisms of En(X)}
holds true for $n=2m+2$, in this case I could not overcome some
technical difficulties which arise in the proof.
\hfill $\bigcirc$
\end{Remark}

\noindent The plan of the proof is as follows. Let $X={\mathbb
C}^m$ or ${\mathbb CP}^m$. To study a strictly equivariant
holomorphic self-map $f$ of the space ${\mathcal E}^n={\mathcal
E}^n(X,gp)$, we start with an explicit description of all
non-constant holomorphic functions $\lambda\colon{\mathcal E}^n
\to \mathbb{C}\setminus \{0,1\}$.\footnote{Compare with
\cite{LinSphere}, \cite{Zinde77a}, \cite{Zinde77b} and
\cite{Feler04a}.} The set $L$ of all such maps is finite and
separates points of a certain rather big submanifold
$M\subset{\mathcal E}^n$ of complex codimension $m(m+1)$. An
endomorphism $f$ induces a self-map $f^*$ of $L$ via
$$
f^*\colon L\ni \lambda
\mapsto f^*\lambda=\lambda\circ f\in L\,.
$$
The map $f^*$ carries important information about $f$.
In order to investigate behaviour of $f^*$
and then recover $f$,
we endow $L$ with the following simplicial structure.
A subset $\Delta^s=\{\lambda_0,...,\lambda_s\}\subseteq L$
is said to be an {\em $s$-simplex} whenever
$\lambda_i/\lambda_j\in L$ for all distinct $i,j$.
The action of ${\mathbf S}(n)$ in ${\mathcal E}^n$
induces a simplicial ${\mathbf S}(n)$ action in
the complex $L$. The orbits of this action may be
exhibited explicitly. On the other hand,
the map $f^*\colon L\to L$ defined above is simplicial and
preserves dimension of simplices. Since $f$ is equivariant,
$f^*$ is nicely related to the ${\mathbf S}(n)$ action on $L$.
\vskip5pt

\noindent Studying all these things together, we find
a holomorphic map
\begin{equation*}
\tau\colon {\mathcal E}^n(X,gp)\to \mathbf{PSL}(m+1,\mathbb{C})
\end{equation*}
and a permutation $\sigma\in{\mathbf S}(n)$
such that $f(q)=\sigma \tau(q)q$.

\begin{Remark}
The topology of the considered spaces is of great interest
by itself. A. I. Barvinok calculated the first homology group of
the ordered space $\Cal{E}^n(\Bbb{C}^2,gp)$, see  \cite{Barvinok86}.
V. Moulton, \cite{Moulton98}, found the generators and some generating
relations of the fundamental groups $\pi_1(\Cal{E}^n(\mathbb{C}^m,gp))$
and $\pi_1(\Cal{E}^n(\mathbb{CP}^m,gp))$.
T. Terasoma, \cite{TerasomaFGM}, proved that this set of
generating relations is complete for the case $m=2$.
\hfill $\bigcirc$
\end{Remark}

\noindent {\bf Acknowledgments:} This paper is based upon my PhD
Thesis supported by Technion. I wish to thank V. Lin who
introduced me to the problem and encouraged me to work on it.
I also wish to express my thanks to M. Entov and A. Vershik for
stimulating discussions and to B. Shapiro for useful references.

\section{Some properties of spaces of generic configurations}
\label{Sec: Some properties of spaces of generic configurations}

\subsection{Notation and definitions}
\label{ss: Notation and definitions}
The spaces $\Cal{E}^n$ of ordered
geometrically generic configurations in $\mathbb{CP}^m$
or ${\mathbb C}^m$ have the following explicit algebraic
description.
\vskip5pt

\noindent Any point $q\in (\mathbb{CP}^m)^n$ may be represented
as a `matrix'
\begin{equation}\label{eq: points as matrices proj}
q=\left(\begin{array}{c}q_1\\ \vdots\\ q_n\end{array}\right)=
\left(\begin{array}{c}\left[ z_{1,1}:\cdots:z_{1,m+1}\right]\\
\vdots\\ \left[z_{n,1}:\cdots: z_{n,m+1}\right]\end{array}
\right)\in (\mathbb{CP}^m)^n\,,
\end{equation}
where $q_j=[z_{j,1}:\cdots:z_{j,m+1}]\in\mathbb{CP}^m$,
$j=1,...,n$.
For $m+1$ distinct indices $i_1,...,i_{m+1}\in \{1,...,n\}$,
the determinant
\begin{equation}\label{def. of 1-cohomology base projective case}
d_{i_1,\dots,i_{m+1}}(q)=
\left|\begin{array}{cccc}  z_{i_1,1}& \dots & z_{i_1,m} & z_{i_1,m+1}
\\ \vdots & \vdots & \vdots &\vdots \\ z_{i_{m+1},1}& \dots &
z_{i_{m+1},m} & z_{i_{m+1},m+1}  \end{array}\right|
\end{equation}
is a homogeneous polynomial of degree $m+1$ in the homogeneous
coordinates
$[z_{1,1}:\cdots:z_{1,m+1}]$,
...,$[ z_{n,1}:\cdots:z_{n,m+1}]$. The space
$\Cal{E}^n({\mathbb{CP}}^m,gp)$ consists of
all matrices $q$ of the form \eqref{eq: points as matrices proj}
such that $d_{i_1,...,i_{m+1}}(q)\ne 0$ for all
distinct $i_1,...,i_{m+1}\in \{1,...,n\}$.
\vskip5pt

\noindent Similarly, the space $\Cal{E}^n({\mathbb C}^m,gp)$
consists of all matrices
\begin{equation}\label{eq: points as matrices aff}
q=\left(\begin{array}{c}q_1\\ \vdots\\ q_n\end{array}\right)
=\left(\begin{array}{c} z_{1,1},...,z_{1,m},1\\
\vdots\\ z_{n,1},...,z_{n,m},1\end{array}\right)
\in (\mathbb C^m)^n
\end{equation}
with all $q_j=(z_{j,1},...,z_{j,m})\in{\mathbb C}^m$
such that
\begin{equation}\label{def. of 1-cohomology base affine case}
d_{i_1,\dots,i_{m+1}}(q)= \left|
\begin{array}{cccc} z_{i_1,1}&
\dots & z_{i_1,m} & 1 \\
\vdots & \vdots & \vdots &\vdots \\
z_{i_{m+1},1}& \dots & z_{i_{m+1},m} & 1
\end{array}
\right|\ne 0
\end{equation}
for all distinct $i_1,...,i_{m+1}\in \{1,...,n\}$.
\vskip5pt

\noindent We refer to the components
$q_1,...,q_n$ of a point
$q=(q_1,...,q_n)$ in $\Cal{E}^n({\mathbb{CP}}^m,gp)$ or
in $\Cal{E}^n({\mathbb C}^m,gp)$
as to {\em vector coordinates} of $q$.
\vskip5pt

\noindent Although we use the same notation for determinant
polynomials in both cases, projective and affine,
every time we meet them it will be clear from the context
which one we mean.
\vskip5pt

\begin{Notation}\label{Not: multiidices}
By a multiindex we mean an ordered set
${\mathbf i}=(i_1,...,i_s)$ with distinct
$i_1,...,i_s\in \{1,...,n\}$.
Sometimes we forget the order and write
$i\in{\mathbf i}$ and $\card {\mathbf i}=s$.
If $s=1$ and $i\in\{1,...,n\}$, we may write
${\mathbf i}=i$ instead of ${\mathbf i}=(i)$.

For $t$ multiindices
${\mathbf i}_1=(i_1^1,...,i_s^1),...,
{\mathbf i}_t=(i_1^t,...,i_s^t)$ such that
$i_1^1,...,i_s^1,...,i_1^t,...,i_s^t$ are distinct,
we set $d_{{\mathbf i}_1,...,{\mathbf i}_t}=
d_{i_1^1,...,i_s^1,...,i_1^t,...,i_s^t}$.

Let ${\mathbf i}=(i_1,...,i_s)$ and ${\mathbf i}'=(i_1,...,i_{s-1})$;
for any $i\in\{i_1,...,i_s\}$ and any $k=1,...,s-1$, we
denote by $D_{{\mathbf i};i,k}$
the $(s-1)\times (s-1)$ minor of the $s\times s$ matrix
$$
Z_{{\mathbf i}}
=\left(\begin{array}{c} z_{i_1,1},...,z_{i_1,s-1},1\\
\vdots\\ z_{i_s,1},...,z_{i_s,s-1},1\end{array}\right)
$$
complementary to the elements $z_{i,k}$; for $i=i_s$,
we write $\delta_{{\mathbf i}';k}$
instead of $D_{{\mathbf i};i_s,k}$.

Let ${\mathbf I}^s$ denote the set of all
multiindices ${\mathbf i}=(i_1,...,i_s)$
such that $1\le i_1<...<i_s\le n$.
For ${\mathbf i}=(i_1,...,i_s)\in {\mathbf I}^s$ and
${\mathbf j}=(j_1,...,j_{t})\in {\mathbf I}^t$,
we define the multiindices ${\mathbf i}\cap{\mathbf j}$ and
${\mathbf i}\setminus {\mathbf j}$ in the evident way;
if ${\mathbf i}\cap{\mathbf j}=\varnothing$,
we define the multiindex
${\mathbf i}\cup{\mathbf j}\in {\mathbf I}^{s+t}$
by an appropriate reordering of the components
$i_1,...,i_s,j_1,...,j_t$.

For $1\le j\le n$, set ${\mathbf I}_j^s=\{{\mathbf i}
\in {\mathbf I}^s\,|\ j \not\in {\mathbf i}\}$.
\hfill $\bigcirc$
\end{Notation}

\subsection{Irreducibility of determinant polynomials}
\label{ss: Irreducibility of determinant polynomials}
\begin{Lemma}\label{cl: irreducibility}
All minors of the matrices
\eqref{eq: points as matrices proj}
and \eqref{eq: points as matrices aff} are irreducible
homogeneous polynomials in the entries $z_{t,s}$.
\end{Lemma}

\begin{proof}
It suffices to prove the lemma for the
matrix \eqref{eq: points as matrices aff}. The proof
is by induction in the order of a minor.
Clearly, any minor of order $1$ is an irreducible polynomial.
Suppose that for some $k$, $1\le k <m+1$,
all minors of order $\le k$ are irreducible.
By the Lagrange decomposition formula,
any minor $M$ of order $k+1$ is a linear function of the
entries $z_{i_1,j_1},...,z_{i_{k+1},j_1}$
of its first column with coefficients that
are certain minors of order $k$, i. e., polynomials of
all entries of $M$ but $z_{i_1,j_1},...,z_{i_{k+1},j_1}$.
By the induction hypothesis, all the latter polynomials
are irreducible; moreover, they depend on different sets of
variables and hence cannot be proportional
(with constant coefficients) to a single polynomial.
This implies that $M$ is irreducible.
\end{proof}

\begin{Lemma}\label{claim reducibility}
Let $i_0\in \{1,...,n\}$, ${\mathbf i}_0=(i_{1},...,i_{m}) \in {\mathbf I}^{m+1}_{i_0}$
and $L \subset \mathbb{C}^{mn}$ be a linear subspace defined by the relations $z_{i_{1},2}=...=z_{i_{m},2}$.

$a)$ For ${\mathbf i}\in {\mathbf I}^{m+1}_{i_0}$, the restriction $d_{{\mathbf i}}|_L$ of $d_{{\mathbf i}}$ to
$L$ is irreducible if and only if
$\card({\mathbf i}\cap{\mathbf i}_0)<m$. Moreover, if $\card({\mathbf i}\cap{\mathbf i}_0)=m$, then
$d_{{\mathbf i}}|_L=\pm D_{{\mathbf i};i,2}\cdot (z_{i,2}-z_{i_1,2})$, where $i\in {\mathbf i}\setminus{\mathbf i}_0$.

$b)$ For ${\mathbf i}\in {\mathbf I}^{m+1}_{i_0}$, $\deg d_{{\mathbf i}}|_L=m$.
\end{Lemma}

\begin{proof}
Assume that $\card({\mathbf i}\cap{\mathbf i}_0)=m$.
It follows that $\card ({\mathbf i}\setminus{\mathbf i}_0)=1$.
Set $(i) = {\mathbf i}\setminus{\mathbf i}_0$.
By the Lagrange determinant decomposition formula,
we can show that
$d_{{\mathbf i}}|_L
=\pm D_{{\mathbf i};i,2}\cdot (z_{i,2}-z_{i_1,2})$.
That is, $d_{{\mathbf i}}|_L$ is reducible and $\deg d_{{\mathbf i}}|_L=m$.
\vskip5pt

\noindent Now assume that
$\card({\mathbf i}\cap{\mathbf i}_0)<m$.
Let $i,j\in {\mathbf i}\setminus{\mathbf i}_0$ be two distinct indices.
It is clear that if we prove that $d_{i,j,i_2,...,i_m}|_L$ is irreducible and
$\deg d_{i,j,i_2,...,i_m}|_L=m$, the same statements for
$d_{{\mathbf i}}|_L$ hold true. So set ${\mathbf j}=(i,j,i_2,...,i_m)$.
Obviously,
$$
d_{{\mathbf j}}|_L
= (z_{i_2,2}-z_{i,2})\cdot D_{{\mathbf j};i,2}-
(z_{i_2,2}-z_{j,2})\cdot D_{{\mathbf j};j,2}
$$
and
$\deg d_{{\mathbf j}}|_L=m$.
It follows that $d_{{\mathbf j}}|_L$ is a linear function of the variables $z_{i,2},z_{j,2}$ and $z_{i_2,2}$ with coefficients
$D_{{\mathbf j};i,2}, D_{{\mathbf j};j,2}$ and
$D_{{\mathbf j};i,2}-D_{{\mathbf j};i,2}$.
By Lemma \ref{cl: irreducibility},
$D_{{\mathbf j};i,2}, D_{{\mathbf j};j,2}$ are irreducible. Therefore
the polynomials
$D_{{\mathbf j};i,2}, D_{{\mathbf j};j,2}$ and
$D_{{\mathbf j};i,2}-D_{{\mathbf j};i,2}$ are pairwise co-prime.
The latter implies that $d_{{\mathbf j}}|_L$ is irreducible.
This completes the proof of the lemma.
\end{proof}

\subsection{The direct decomposition of $\Cal{E}^n(\Bbb{CP}^m,gp)$}
\label{ss: The direct decomposition of En(CPm)}
Here we show that $\Cal{E}^n(\Bbb{CP}^m,gp)$ admits
a natural representation as a Cartesian product of
its subspace of codimension $m(m+2)$
and the group $\mathbf{PSL}(m+1,\mathbb{C})$.

\begin{Definition}\label{Def: M subspace}
Set
\begin{equation}\label{eq: vectors v}
\aligned
&v_1=[1:0:\cdots:0]\,, \ v_2=[0:1:0:\cdots:0]\,,...,
v_{m+1}=[0:\cdots:0:1] \\
&\text{and} \ \ w=[1:\cdots:1]\,.
\endaligned
\end{equation}
The subspace $M_{m,n}\subset\Cal{E}^n(\Bbb{CP}^m,gp)$
defined by
$$
M_{m,n} = \left\{q=(q_1,...,q_n)\in
 \Cal{E}^n(\Bbb{CP}^m,gp)|\
 q_i = v_i \ \forall \, i=1,...,m+1\,,\
 q_{m+2}=w
 \right\}
$$
is called the {\em reduced space
of geometrically generic ordered configurations}.
\hfill $\bigcirc$
\end{Definition}

\begin{Lemma}
\label{Lm: existence of the special map in the space}
Let $n\ge m+3$. For every $q\in\Cal{E}^n(\mathbb{CP}^m,gp)$,
there is a unique $\gamma(q)\in\mathbf{PSL}(m+1,\mathbb{C})$ such that
$\gamma(q)q \in M_{m,n}$. The map
$$
\gamma\colon\Cal{E}^n(\Bbb{CP}^m,gp)\ni q
\mapsto\gamma(q)\in\mathbf{PSL}(m+1,\mathbb{C})
$$
is holomorphic.
\end{Lemma}

\begin{proof}\footnote{\noindent For $m=1$ the statement of Lemma
is a common knowledge; the case $m=2$ is treated
in \cite{Efimov(rus)}, Chap. V, Sec. 109, Theorem 36.
For $m>2$, I could not find an
appropriate reference and sketched the proof here.}
For
$
q=\left(\begin{array}{c}\left[ z_{1,1}:\cdots:z_{1,m+1}\right]\\
\vdots\\ \left[z_{n,1}:\cdots: z_{n,m+1}\right]\end{array}
\right)
\in \Cal{E}^n(\Bbb{CP}^m,gp)
$ and $i\in\{1,...,m+1\}$,
the matrices
$$
A(q)=\left(\begin{array}{lcr}z_{1,1}&,...,&z_{1,m+1}\\
&\vdots&\\ z_{m+1,1}&,...,& z_{m+1,m+1}\end{array}
\right)\,,  \ A_i(q)=\left(\begin{array}{lcr}z_{1,1}&,
...,&z_{1,m+1}\\
&\vdots&\\
z_{i-1,1}&,...,&z_{i-1,m+1}\\
z_{m+2,1}&,...,&z_{m+2,m+1}\\
z_{i+1,1}&,...,&z_{i+1,m+1}\\
&\vdots&\\
z_{m+1,1}&,...,&z_{m+1,m+1}\end{array}
\right)
$$
are defined up to multiplication of their rows by non-zero
complex numbers. The determinants $\det A(q)$ and all $D_i=\det A_i$
are homogeneous polynomials non-vanishing on
$\Cal{E}^n(\Bbb{CP}^m,gp)$; therefore, the determinant of the
adjunct matrix $X(q)=A^{adj}(q)$
(whose elements $x_{i,j}$ are the algebraic co-factors of
the elements $z_{j,i}$ of $A(q)$) is equal to $(\det A(q))^m$ and
hence does not vanish on $\Cal{E}^n(\Bbb{CP}^m,gp)$.
It follows that for every $q\in\Cal{E}^n(\Bbb{CP}^m,gp)$
the matrix $T(q)=(x_{i,j}/D_j)_{i,j=1}^{m+1}$
is non-degenerate. It is easy to verify that $T(q)$
determines a well-defined automorphism $\gamma(q)$ of $\mathbb{CP}^m$
and $\gamma(q)q\in M_{m,n}$.

An automorphism of $\Bbb{CP}^m$ is uniquely
determined by its values at any $m+2$ points in geometrically
generic position. It follows that for any
$q\in\Cal{E}^n(\Bbb{CP}^m,gp)$ the constructed above element
$\gamma(q)\in\mathbf{PSL}(m+1,\mathbb{C})$ is the only element of
$\mathbf{PSL}(m+1,\mathbb{C})$ that carries $q$ to a point of $M_{m,n}$.
Moreover, $\gamma(q)$ holomorphically depends on a point $q$.
\end{proof}

\begin{Corollary}\label{Crl: direct decomposition for projective case}
The mutually inverse maps
$$
A\colon\Cal{E}^n(\Bbb{CP}^m,gp)\ni q \mapsto A(q)=
(\gamma(q),\gamma(q)q)\in \mathbf{PSL}(m+1,\mathbb{C})\times M_{m,n}
$$
and
$$
B\colon\mathbf{PSL}(m+1,\mathbb{C})\times
M_{m,n}\ni(T,\tilde{q})\mapsto
B(T,\tilde{q})= T^{-1}\tilde{q}\in
\Cal{E}^n(\Bbb{CP}^m,gp)
$$
induce a natural biholomorphic isomorphism
$\Cal{E}^n(\Bbb{CP}^m,gp)\cong\mathbf{PSL}(m+1,\mathbb{C})\times M_{m,n}$.
\end{Corollary}

\begin{Remark}\label{rm: configuration spaces are affine}
The above corollary implies that $\Cal{E}^n(\Bbb{CP}^m,gp)$
and $\Cal{C}^n(\Bbb{CP}^m,gp)$ are irreducible non-singular
complex affine algebraic varieties and, in particular,
Stein manifolds. Indeed, in the above direct decomposition
of $\Cal{E}^n(\Bbb{CP}^m,gp)$ both $\mathbf{PSL}(m+1,\mathbb{C})$ and
$M_{m,n}$ are such varieties and hence $\Cal{E}^n(\Bbb{CP}^m,gp)$
is also of the same nature. Since the group ${\mathbf S}(n)$
is finite and its action on $\Cal{E}^n(\Bbb{CP}^m,gp)$ is free,
the same is true for the quotient $\Cal{C}^n(\Bbb{CP}^m,gp)
=\Cal{E}^n(\Bbb{CP}^m,gp)/{\mathbf S}(n)$.

The same properties hold true for $\Cal{E}^n(\Bbb{C}^m,gp)$
and $\Cal{C}^n(\Bbb{C}^m,gp)$.

Notice that for $n>1$ the ``standard" $n$-point configuration
spaces $\Cal{E}^n(X)$ and $\Cal{C}^n(X)$ of a complex manifold $X$
of dimension $m>1$ cannot be Stein manifolds.
\hfill $\bigcirc$
\end{Remark}

\subsection{Determinant cross ratios}
\label{ss: Determinant cross ratios}
Here we construct certain non-constant holomorphic functions
$\Cal{E}^n\to{\mathbb C}\setminus\{0,1\}$, which are called
``determinant cross ratios" (in fact, later we shall show
that there are no other functions with these properties).

Let us recall that, according to
Notation \ref{Not: multiidices},
for distinct $i_1,...,i_{m-1},j,k$ and
${\mathbf i}=(i_1,...,i_{m-1})$
the notation $d_{{\mathbf i},j,k}$ means the determinant
$d_{i_1,...,i_{m-1},j,k}$.

\begin{Definition}\label{Def: determinant cross ratio}
Let $X$ be either $\mathbb{CP}^m$ or $\mathbb{C}^m$ and let
$n\ge m+3$. For any $m+3$-dimensional multiindex
$I=(i_1,...,i_{m+3})$ with distinct components $i_t\in\{1,...,n\}$, set
${\mathbf i}=(i_1,...,i_{m-1})$,
$j=i_m$, $k=i_{m+1}$, $l=i_{m+2}$, $s=i_{m+3}$;
the non-constant rational function
\begin{equation}\label{eq: def of e}
e_I(q)=e_{{\mathbf i};j,k,l,s}(q)=
\frac{d_{{\mathbf i},j,k}(q)}{d_{{\mathbf i},j,l}(q)}:
\frac{d_{{\mathbf i},k,s}(q)}{d_{{\mathbf i},l,s}(q)}\,, \ \
q\in ({\mathbb{CP}}^m)^n\,,
\end{equation}
is called a {\em determinant cross ratio}, or, in brief,
a DCR. This function is regular on the algebraic manifold
$\Cal{E}^n(X,gp)\subset ({\mathbb{CP}}^m)^n$.

The unordered set of indices
$\{I\}=\{i_1,\dots,i_{m-1},j,k,l,s\}$ is called the
{\em support} of the function $\mu=e_I=e_{{\mathbf i};j,k,l,s}$
and is denoted by $\supp\mu$; its unordered subset
$\{\mathbf i\}=\{i_1,...,i_{m-1}\}$
is called the {\em essential support} of $\mu$
and is denoted by $\supp_{\text{ess}}\mu$ (we often write
$I$ instead of $\{I\}$ and $\mathbf i$ instead of $\{\mathbf i\}$).
In fact, the function $e_I(q)=e_I(q_1,...,q_n)$
depends only on the vector variables $q_t$ with $t\in I$.
\hfill $\bigcirc$
\end{Definition}

\begin{Remark}\label{Rmk: when two DCR coincide}
Notice that two determinant cross ratios, say
$e_I=e_{{\mathbf i};j,k,l,s}$
and $e_{I'}=e_{{\mathbf i}';j',k',l',s'}$,
coincide if and only if $\{{\mathbf i}\}=\{{\mathbf i}'\}$
and $(j',k',l',s')$ is obtained from $(j,k,l,s)$ by a
Kleinian permutation of four letters.
The set of all determinant cross ratios
is denoted by $\dcr (\Cal{E}^n)=\dcr (\Cal{E}^n(X,gp))$.
\hfill $\bigcirc$
\end{Remark}

\noindent The ${\mathbf S}(n)$ action in ${\mathcal E}^n$
induces an ${\mathbf S}(n)$ action on functions defined by
$(\sigma\lambda)(q)=\lambda(\sigma^{-1}q)=
\lambda(q_{\sigma(1)},...,q_{\sigma(n)})$ \
($\lambda$ is a function on ${\Cal E}^n$,
$q=(q_1,...,q_n)\in{\Cal E}^n$; notice that
$\sigma\lambda$ may also be written as
$(\sigma^{-1})^*\lambda=\lambda\circ(\sigma^{-1})$,
where $\sigma$ and $\sigma^{-1}$ are considered as the
self-mappings of ${\mathcal E}^n$). This action
carries holomorphic functions to holomorphic functions.

\begin{Lemma}\label{Rmk: transitivity of cross ratios}
For any $\sigma\in \mathbf{S}(n)$ and any
$\lambda\in\dcr (\Cal{E}^n)$ the function
$\sigma\lambda$ also belongs to $\dcr (\Cal{E}^n)$.
Moreover, the ${\mathbf S}(n)$ action is transitive on the set
$\dcr (\Cal{E}^n)$.
\end{Lemma}

\begin{proof}
Let $\mathbf{a}=(a_1,...,a_{l})$ be a multiindex.
For any $\sigma \in {\mathbf S}(n)$, let $\sigma (\mathbf{a})$
denote the multiindex $(\sigma  (a_1),...,\sigma  (a_{l}))$.
For any $\sigma \in {\mathbf S}(n)$, we have
$\sigma d_{\mathbf{a}}(q)=d_{\mathbf{a}}(\sigma^{-1}q)
=d_{\sigma(\mathbf{a})}(q)$;
thus,
$$
\aligned
\sigma e_{{{\mathbf i}};j,k,r,s}(q)&=
e_{{{\mathbf i}};j,k,r,s}(\sigma^{-1} q)=
    \frac{d_{{\mathbf i},j,k}(\sigma^{-1} q)}
    {d_{{\mathbf i},j,r}(\sigma^{-1} q)}:
  \frac{d_{{\mathbf i},k,s}(\sigma^{-1} q)}
  {d_{{\mathbf i},r,s}(\sigma^{-1} q)}\\
  &=
  \frac{d_{\sigma({\mathbf i}),\sigma(j),\sigma(k)}(q)}
  {d_{\sigma({\mathbf i}),\sigma(j),\sigma(r)}(q)}:
  \frac{d_{\sigma({\mathbf i}),\sigma(k),\sigma(s)}(q)}
  {d_{\sigma({\mathbf i}),\sigma(r),\sigma(s)}(q)}=
  e_{{\sigma({\mathbf i})};\sigma(j),
  \sigma(k),\sigma(r),\sigma(s)}(q)\,.
  \endaligned
  $$
Let $e_I,e_{I'}\in\dcr({\mathcal E}^n)$.
Since each of the ordered sets $I$ and $I'$ consists of
$m+3\le n$ distinct elements of the set $\{1,...,n\}$,
there is a permutation $\sigma\in{\mathbf S}(n)$
such that $\sigma I= I'$ and hence $\sigma e_I=e_{I'}$.
\end{proof}

\begin{Lemma}\label{Lm: invariants of psl(m+1) action}
Determinant cross ratios are invariants of the
$\mathbf{PSL}(m+1,\mathbb{C})$ action on the configuration space
$\Cal{E}^n(\Bbb{CP}^m,gp)$.
\end{Lemma}

\begin{proof}
It is easy to observe that the following elementary
operators do not change any determinant cross ratio:
\begin{align*}
&\left[ z_{1} : \dots : z_{m+1}\right]
\mapsto
\left[ a_{1} z_{1}: \dots : a_{m+1} z_{m+1}\right]
\quad  \text{ for } a_{1}\cdot...\cdot a_{m+1}\not=0\,;
\\
&\left[z_{1}:\dots :z_{i}:\dots:z_{j}:\dots:z_{m+1}\right]
\mapsto
\left[z_{1}:\dots: z_{j}:\dots:z_{i}:\dots:z_{m+1}\right]\,;
\\
&\left[z_{1}:\dots:z_{i}:\dots:z_{j}:\dots:z_{m+1}\right]
\mapsto
\left[z_{1}:\dots:z_{i}+z_{j}:\dots:z_{j}:\dots:z_{m+1}\right]\,.
\end{align*}
Any element of $\mathbf{PSL}(m+1,\mathbb{C})$ can be decomposed into
a sequence of the elementary operators. This proves the lemma.
\end{proof}

\begin{Lemma}\label{Lm: determinant cross ratios omit 0 and 1}
Any function $\lambda \in \dcr({\mathcal E}^n)$
omits the values $0$ and $1$.
\end{Lemma}

\begin{proof}
By Lemma \ref{Lm: invariants of psl(m+1) action},
any cross ratio $\lambda$ is $\mathbf{PSL}(m+1,\mathbb{C})$ invariant.
By Lemma \ref{Lm: existence of the special map in the space},
any orbit of $\mathbf{PSL}(m+1,\mathbb{C})$ action in
$\Cal{E}^n$ intersects the subspace $M_{m,n}$;
hence, it sufices to show that $\lambda$ does not accept
the values $0$ and $1$ on $M_{m,n}$.

First, let
$I=(1,....,m+3)$, ${\mathbf i}=(1,...,m-1)$ and $\lambda=e_I(q)$;
then
$$
\lambda(q)=e_I(q)
=\frac{d_{{\mathbf i},m,m+1}(q)}
{d_{{\mathbf i},m,m+2}(q)}:
\frac{d_{{\mathbf i},m+1,m+3}(q)}
{d_{{\mathbf i},m+2,m+3}(q)}\,.
$$
For $q=(q_1,...,q_n)\in M_{m,n}$ with
$q_t=[z_{t,1}:\cdots:z_{t,m+1}]$, we have
$$
\aligned
d_{{\mathbf i},m,m+1}(q)&=1\,, \hskip64.5pt
d_{{\mathbf i},m,m+2}(q)=1\,,\\
d_{{\mathbf i},m+1,m+3}(q)&= - z_{m+3,m}\,, \hskip20pt
d_{{\mathbf i},m+2,m+3}(q) = z_{m+3,m+1}-z_{m+3,m}\,,
\endaligned
$$
and hence
$$
e_I(q)=1-\frac{z_{m+3,m+1} }{ z_{m+3,m}}\,.
$$
Since $d_{\mathbf{i},m,m+3}(q)\ne 0$,
$d_{\mathbf{i},m+1,m+3}(q)\ne 0$ and
$d_{\mathbf{i},m+2,m+3}(q)\ne 0$,
we have\linebreak[4] $z_{m+3,m},z_{m+3,m+1}\ne 0$ and
$z_{m+3,m}\ne z_{m+3,m+1}$. Thus, $e_I(q)\ne 0,1$
on $M_{m,n}$ and, consequently, on the whole of ${\mathcal E}^n$.

By Lemma \ref{Rmk: transitivity of cross ratios},
${\mathbf S}(n)$ acts transitively on the set
$\dcr({\mathcal E}^n)$, which implies that any
$\lambda\in\dcr({\mathcal E}^n)$ does not accept the values
$0$ and $1$.
\end{proof}

\begin{Notation}\label{Nt: first m members}
For $s\in \{1,...,m\}$, set
${\mathbf m}(\hat s) = (1,...,\hat{s}, ...,m)$.
For $s=m$, we write sometimes $\widehat{\mathbf m}$
instead of ${\mathbf m}(\hat m)$.
\end{Notation}


\noindent We have also the following lemma:

\begin{Lemma}\label{Lm: definition and properties of P}
$a)$ The map $P\colon M_{m,n}\to (\Bbb{C}^{n-m-2})^m$ defined by
\begin{equation*}
q \mapsto P(q)=\begin{pmatrix}
&p_{1,m+3}(q),&...,&p_{1,n}(q)\\
&... &... &...\\
&p_{m,m+3}(q),&...,&p_{m,n}(q)\\
\end{pmatrix}\,,
\end{equation*}
with $p_{s,t}(q)= e_{{\mathbf m}(\hat s);s,m+1,m+2,t}(q)$
for $s=1,...,m$ and $t=m+3,...,n$,
is a holomorphic embedding.\\
$b)$ $M_{m,n}$ is a hyperbolic space.
\end{Lemma}

\begin{proof}
$a)$ For $q=(q_1,...,q_n)\in M_{m,n}$ with
$q_i=[z_{i,1}:\cdots:z_{i,m+1}]$
we have $z_{t,m+1}=d_{\widehat{\mathbf m},m,t}(q)\ne 0$
and $z_{t,s} =\pm  d_{{\mathbf m}(\hat s),m+1,t}(q)\ne 0$.
Furthermore,
$$
\aligned
d_{{\mathbf m}(\hat s),s,m+1}(q)&=(-1)^{m-s}\,, \hskip44.5pt
d_{{\mathbf m}(\hat s),s,m+2}(q)=(-1)^{m-s}\,,\\
d_{{\mathbf m}(\hat s),m+1,t}(q)&= -(-1)^{m-s} z_{t,m}\,,\hskip20pt
d_{{\mathbf m}(\hat s),m+2,t}(q)
= (-1)^{m-s}(z_{t,m+1}-z_{t,m})\,;
\endaligned
$$
thus,
$$
p_{s,t}(q)=e_{{\mathbf m}(\hat s);s,m+1,m+2,t}(q) =
\frac{d_{{\mathbf m}(\hat s),s,m+1}(q)}
{d_{{\mathbf m}(\hat s),s,m+2}(q)}:
\frac{d_{{\mathbf m}(\hat s),m+1,t}(q)}
{d_{{\mathbf m}(\hat s),m+2,t}(q)}
=  1-\frac{z_{t,m+1} }{ z_{t,s}}\,.
$$
If $q'=(q'_1,...,q'_n)\in M_{m,n}$ with
$q'_i=[z'_{i,1}:\cdots:z'_{i,m+1}]$
and $p_{s,t}(q)=p_{s,t}(q')$ for $s=1,...,m$, then
$z_{t,m+1}:z_{t,s}=z'_{t,m+1}:z'_{t,s}$
for all $s=1,...,m$ and hence $q_t=q'_t$. Thus,
$P(q)=P(q')$ implies $q=q'$ and $P$ is injective.

To see that $P$ is an embedding, it suffices to observe that
from the above calculation follows that at any point
the Jacobi matrix of $P$ is of maximal rank.\\
$b)$ Since every determinant cross ratio omits values $0$ and $1$,
$$
P(M_{m,n})\subset ({\mathbb C}\setminus \{0,1\})^{m(n-m-2)}\,,
$$
it follows that $M_{m,n}$ is a Kobayashi's hyperbolic complex manifold.
\end{proof}

\noindent Lemma \ref{Lm: definition and properties of P}
and Corollary \ref{Crl: direct decomposition for projective case}
imply that the determinant cross ratios generate the whole
algebra ${\mathcal A}
={\mathbb C}[\Cal{E}^n(\Bbb{CP}^m,gp)]^{\mathbf{PSL}(m+1,\mathbb{C})}
={\mathbb C}[M_{m,n}]$
of invariants of the $\mathbf{PSL}(m+1,\mathbb{C})$ action on the
algebraic variety $\Cal{E}^n(\Bbb{CP}^m,gp)$.
We shall see that the set $\dcr$ of all these generators
coincides with the set $L({\mathcal E}^n)$ of all non-constant
holomorphic functions
${\mathcal E}^n\to{\mathbb C}\setminus\{0,1\}$ and that
a strictly equivariant holomorphic
map $f\colon\Cal{E}^n(\Bbb{CP}^m,gp)\to\Cal{E}^n(\Bbb{CP}^m,gp)$
induces a selfmap $f^*$ of $L({\mathcal E}^n)$;
thereby, such a map $f$ induces an endomorphism of the algebra
${\mathcal A}$. Together with Corollary
\ref{Crl: direct decomposition for projective case},
this provides an important information about $f$ and
eventually leads to the proof of Theorem
\ref{Thm: endomorphisms of En(X)}.

\section{Holomorphic functions omitting the values $0$ and $1$}
\label{Sec: Holomorphic functions omitting two values}

\noindent In this section we describe explicitly all holomorphic
functions $\mu\colon\Cal{E}^n\to{\mathbb C}\setminus\{0,1\}$.
\vskip5pt

\subsection{ABC lemma}\label{ss: ABC lemma}
The following lemma plays a crucial part
in an explicit description of all holomorphic
functions $\Cal{E}^n(\Bbb{C}^m,gp)
\to{\mathbb C}\setminus\{0,1\}$.
Let us recall that, according to
Notation \ref{Not: multiidices},
for distinct $i_1,...,i_{m-1},j,k$ and ${\mathbf i}=(i_1,...,i_{m-1})$
the notation $d_{{\mathbf i},j,k}$ means the determinant
$d_{i_1,...,i_{m-1},j,k}$.

\begin{Lemma}\label{Lm: likeABC for greater dimensions}
Let $ A,B,C \in {\mathbb C}[\Bbb{C}^{mn}]
={\mathbb C}[z_{1,1},...,z_{n,m}]$
be pairwise co-prime polynomials on $\Bbb{C}^{mn}$
non-vanishing on the configuration space
$\Cal{E}^n(\Bbb{C}^m,gp)\subset \Bbb{C}^{mn}$.
Assume that at least one
of them is non-constant and $A+B+C=0$. Then there exist
a multiindex $\mathbf{i}=(i_1,\dots,i_{m-1})$, indices $j,k,l,s$,
and a complex number $\alpha\ne 0$ such that
all $i_1,\dots,i_{m-1},j,k,l,s$ are distinct and
$A=\alpha d_{\mathbf{i},j,k} d_{\mathbf{i},l,s}$,
$B=\alpha d_{\mathbf{i},j,l} d_{\mathbf{i},s,k}$ and
$C=\alpha d_{\mathbf{i},j,s} d_{\mathbf{i},k,l}$.
\end{Lemma}

\begin{proof}
The polynomials $A,B,C$ do not vanish on
$$
\Cal{E}^n(\Bbb{C}^m,gp)=\Bbb{C}^{mn}\setminus
\bigcup\limits_{\mathbf i}\{q\in\Bbb{C}^{mn}\,|
\ d_{\mathbf i}(q)=0 \}\,.
$$
It follows from Lemma \ref{cl: irreducibility} that
\begin{gather*}
A=\alpha \prod\limits_{{\mathbf i}\in {\mathbf I}^{m+1}} d_{{\mathbf i}}^{a_{{\mathbf i}}}\,, \quad
B=\beta \prod\limits_{{\mathbf i}\in {\mathbf I}^{m+1}} d_{{\mathbf i}}^{b_{{\mathbf i}}}\,, \quad
C=\gamma \prod\limits_{{\mathbf i}\in {\mathbf I}^{m+1}} d_{{\mathbf i}}^{c_{{\mathbf i}}}\,,
\end{gather*}
where $\alpha,\beta,\gamma \in \mathbb{C} \setminus \{0\}$ and
$a_{{\mathbf i}},b_{{\mathbf i}},c_{{\mathbf i}}\in{\mathbb Z_+}$.
The polynomials are homogeneous; thus, the equality $A+B+C=0$
implies that $\deg A = \deg B =\deg C$; i.e.,
$\sum a_{{\mathbf i}}=\sum b_{{\mathbf i}}=\sum c_{{\mathbf i}}$.
For every index $i_0$, $1\le i_0 \le n$, we can write
\begin{gather*}
A=A_{{i_0}}\prod\limits_{{\mathbf i}\in {\mathbf I}^{m}_{i_0}} d_{{i_0},{\mathbf i}}^{a_{{i_0},{\mathbf i}}}\,,\qquad
B= B_{{i_0}}\prod\limits_{{\mathbf i}\in {\mathbf I}^{m}_{i_0}} d_{{i_0},{\mathbf i}}^{b_{{i_0},{\mathbf i}}}\,,\qquad
C=C_{{i_0}}\prod\limits_{{\mathbf i}\in {\mathbf I}^{m}_{i_0}} d_{{i_0},{\mathbf i}}^{c_{{i_0},{\mathbf i}}}\,,
\end{gather*}
where $A_{{i_0}},B_{{i_0}},B_{{i_0}}$ are the products of all
factors $d_{\mathbf i}$ that do not contain the variables
$z_{i_0,1},...,z_{i_0,m}$, i. e.,
\begin{equation}\label{eq: present of Ai0 Bi0 Ci0}
A_{{i_0}}=\pm\alpha \prod\limits_{{\mathbf i}\in {\mathbf I}^{m+1}_{i_0}}
d_{{\mathbf i}}^{a_{{\mathbf i}}}\,,\quad
B_{{i_0}}=\pm\beta \prod\limits_{{\mathbf i}\in {\mathbf I}^{m+1}_{i_0}}
d_{{\mathbf i}}^{b_{{\mathbf i}}}\,,\quad
C_{{i_0}}=\pm\gamma \prod\limits_{{\mathbf i}\in {\mathbf I}^{m+1}_{i_0}}
d_{{\mathbf i}}^{c_{{\mathbf i}}}\,.
\end{equation}
Our main tool is induction on the dimension $m$.
The proof is divided into two steps;
the first one supplies us with
a way of the reduction of the dimension,
and the second one is the induction itself.
\vskip7pt

\noindent {\bf Step 1.} Let us prove the following statement:
\vskip5pt

$(*)$ \emph{There is an index $t_0$ such that the polynomials
$A_{t_0},B_{t_0},C_{t_0}$ are constant.}
\vskip5pt

{\it Proof.}
First, we prove that for any $i_0\in \{1,...,n\}$
\begin{equation}\label{claim 4}
\sum a_{{i_0},{\mathbf i}}=
\sum  b_{{i_0},{\mathbf i}} = \sum c_{{i_0},{\mathbf i}} \
\text{and} \  \deg A_{{i_0}} = \deg B_{{i_0}} =\deg C_{{i_0}}\,.
\end{equation}
Without loss of generality, we can
assume that either $a)$ $\sum a_{{i_0},{\mathbf i}}>
\sum b_{{i_0},{\mathbf i}}\ge\sum c_{{i_0},{\mathbf i}}$
or $b)\, \sum a_{{i_0},{\mathbf i}}=\sum b_{{i_0},{\mathbf i}}>
\sum c_{{i_0},{\mathbf i}}$ or
$c)\, \sum a_{{i_0},{\mathbf i}}=\sum b_{{i_0},{\mathbf i}}=
\sum c_{{i_0},{\mathbf i}}$.

Compare the terms of the maximal degree in the variable
$ z_{{i_0},1}$ in the main equality $A+B+C=0$.
In the case $(a)$ we have
$\displaystyle A_{{i_0}}\prod \delta_{{\mathbf i};1}^{a_{i_0,{\mathbf i}}}=0$
(see Notation \ref{Not: multiidices}
for the definition of $\delta_{{\mathbf i};k}$).
This means $A=0$, a contradiction.
In the case $(b)$ we obtain
$\displaystyle
A_{{i_0}}\prod\limits_{{\mathbf i}\in {\mathbf I}^{m}_{i_0}} \delta_{{\mathbf i};1}^{a_{{i_0},{\mathbf i}}} +B_{{i_0}}\prod\limits_{{\mathbf i}\in {\mathbf I}^{m}_{i_0}} \delta_{{\mathbf i};1}^{b_{{i_0},{\mathbf i}}}=0.
$
By Lemma \ref{cl: irreducibility}, we have\footnote{We denote the greatest common divisor by $\Gcd$.}
$$
\Gcd(A_{{i_0}},\prod\limits_{{\mathbf i}\in {\mathbf I}^{m}_{i_0}}
\delta_{{\mathbf i};1}^{b_{{i_0},{\mathbf i}}})=\Gcd(B_{{i_0}},
\prod\limits_{{\mathbf i}\in {\mathbf I}^{m}_{i_0}}
\delta_{{\mathbf i};1}^{a_{{i_0},{\mathbf i}}})=1\,.
$$
It follows that
$A_{{i_0}}+B_{{i_0}}=0$ and $a_{{i_0},{\mathbf i}}=b_{{i_0},
{\mathbf i}}$ for all ${\mathbf i}\in {\mathbf I}^{m}_{i_0}$.
The latter implies that $C=-(A+B)=0$,
a contradiction. This completes the proof of \eqref{claim 4}.
\vskip5pt

\noindent Since $A+B+C=0$, it follows that for each $k=1,...,m$
the leading term of $A+B+C$ in the variable $z_{{i_0},k}$ is $0$.
In other words, for each $k=1,...,m$, we obtain the following equality:
\begin{equation}\label{eq:needs Lin's theorem1 general}
A_{{i_0}}\prod\limits_{{\mathbf i}\in {\mathbf I}^{m}_{i_0}} \delta_{{\mathbf i};k}^{a_{{i_0},{\mathbf i}}} +B_{{i_0}}\prod\limits_{{\mathbf i}\in {\mathbf I}^{m}_{i_0}} \delta_{{\mathbf i};k}^{b_{{i_0},{\mathbf i}}}+
C_{{i_0}}\prod\limits_{{\mathbf i}\in {\mathbf I}^{m}_{i_0}} \delta_{{\mathbf i};k}^{c_{{i_0},{\mathbf i}}}=0\,.
\end{equation}

We may assume that there is an index ${i_0}$ such that
$A_{i_0}\ne A$ and $A_{i_0}\ne \const$; for otherwise,
either $A=\const$, which is a contradiction, or $A_{i_0}=\const$
and $(*)$ holds with $t_0=i_0$.
According to \eqref{claim 4}, for such $i_0$ we have
$0<\deg A_{i_0}=\deg B_{i_0}=\deg C_{i_0} < \deg A=\deg B=\deg C$.

By Lemma \ref{cl: irreducibility},
the polynomials $\delta_{{\mathbf i};k}$ are irreducible and
distinct. Since $A,B,C$ are pairwise co-prime, the integers
$a_{i_0,{\mathbf i}},b_{i_0,{\mathbf i}},c_{i_0,{\mathbf i}}$
are also distinct.
Using these facts, it easy to verify that
the second order system of the linear equations in variables
$A_{i_0}, B_{i_0},C_{i_0}$ defined by
\eqref{eq:needs Lin's theorem1 general} with $k=1,2$ is
of rank $2$; actually, all its $2\times 2$ minors are
non-zero polynomials. It follows that
$\displaystyle\frac{A_{i_0}}{\mathfrak{A}}=
\frac{B_{i_0}}{\mathfrak{B}}=\frac{C_{i_0}}{\mathfrak{C}}$
with certain non-zero polynomials
$\mathfrak{A},\mathfrak{B},\mathfrak{C}$.
The polynomials $B_{i_0}$ and $C_{i_0}$ are co-prime, consequently,
\begin{equation}\label{eq:case3 general}
\mathfrak{B}= \tilde{B}\cdot B_{i_0} \ \ \text{and} \ \
\mathfrak{C}= \tilde{C}\cdot C_{i_0}\,,
\end{equation}
where $\tilde{B},\tilde{C}$ are non-zero polynomials.
Pick a multiindex ${\mathbf i}_0=(i_1,\dots,i_m)\in
{\mathbf I}^{m}_{i_0}$ such that $ a_{{i_0},{\mathbf i}_0}>0$.
Since the polynomials $A,B,C$ are pairwise co-prime,
it follows that
$ b_{{i_0},{\mathbf i}_0}=c_{{i_0},{\mathbf i}_0}=0$.
\vskip5pt

\noindent Let us prove that
\begin{equation}\label{eq:for A 1 general}
A = A_{{i_0}} d_{{i_0},{\mathbf i}_0}^{a_{{i_0},{\mathbf i}_0}}\,,\qquad
B_{{i_0}} = \pm \beta \prod\limits_{i\not={i_0}} d_{i,{\mathbf i}_0}^{b_{i,{\mathbf i}_0}}\,,\qquad
C_{{i_0}} = \pm\gamma\prod\limits_{i\not={i_0}} d_{i,{\mathbf i}_0}^{c_{i,{\mathbf i}_0}}\,,
\end{equation}
the statement $(*)$ will follow from this by combinatorial
considerations.

Let $L$ be the linear subspace of $\mathbb{C}^{mn}$
defined by the relations $z_{i_{1},2}=...=z_{i_{m},2}$.
Since for ${\mathbf i}\in {\mathbf I}^{m}_{i_0}$ we have
$\delta_{{\mathbf i};k}|_L=0$ if and only if $\card({\mathbf i}\cap{\mathbf i}_0)=m$ and $k\ne2$,
the restrictions
\begin{equation}\label{eq: restricted B and C}
\mathfrak{B}|_L= -\left.\left(\prod
\delta_{{\mathbf i};2}^{a_{{i_0},{\mathbf i}}}
\cdot\prod
\delta_{{\mathbf i};1}^{c_{{i_0},{\mathbf i}}}\right)\right|_L \quad
\text{\rm and}
\quad
\mathfrak{C}|_L= -\left.\left(\prod
\delta_{{\mathbf i};2}^{a_{{i_0},{\mathbf i}}}\cdot
\prod
\delta_{{\mathbf i};1}^{b_{{i_0},{\mathbf i}}}\right)\right|_L
\end{equation}
are non-zero polynomials.
Thus $\tilde{B}|_L,\tilde{C}|_L\ne 0$.
According to \eqref{eq: restricted B and C},
the polynomials $\mathfrak{B}|_L$ and
$\mathfrak{C}|_L$ are products of irreducible polynomials
of degree $\le (m-1)$. By \eqref{eq: present of Ai0 Bi0 Ci0},
\begin{equation}\label{eq: B0 C0 decompositions on L}
B_{i_0}|_L=\pm\beta \prod\limits_{{\mathbf i}\in
{\mathbf I}^{m+1}_{i_0}}(d_{{\mathbf i}}|_L)^{b_{{\mathbf i}}}
\qquad \text{and} \qquad
C_{i_0}|_L=\pm\gamma \prod\limits_{{\mathbf i}\in {\mathbf I}^{m+1}_{i_0}}
(d_{{\mathbf i}}|_L)^{c_{{\mathbf i}}}\,.
\end{equation}
By \eqref{eq:case3 general}, $B_{i_0}|_L$ and $C_{i_0}|_L$
also must be products of irreducible polynomials
of degree $\le (m-1)$.
By Lemma \ref{claim reducibility}, this may happen if and only if
the decompositions
\eqref{eq: B0 C0 decompositions on L}
of $B_{i_0}|_L$ and $C_{i_0}|_L$
contain only factors of the form $ d_{i,{\mathbf i}_0}|_L$.
I.e., we have showed that
$$
B_{{i_0}} = \pm \beta \prod\limits_{i\not={i_0}}
d_{i,{\mathbf i}_0}^{b_{i,{\mathbf i}_0}}\quad
\text{and} \quad
C_{{i_0}} = \pm\gamma\prod\limits_{i\not={i_0}}
d_{i,{\mathbf i}_0}^{c_{i,{\mathbf i}_0}}\,.
$$
To complete the proof of
\eqref{eq:for A 1 general}, we must show that
$a_{i_0,{\mathbf i}_1}=0$ for any ${\mathbf i}_1\in
{\mathbf I}^m_{i_0}\setminus\{{\mathbf i}_0\}$.
Suppose on the contrary that there exists
a multiindex ${\mathbf i}_1$ such that
$a_{{i_0},{\mathbf i}_1}>0$ and
${\mathbf i}_0\not={\mathbf i}_1$.
Then the same argument implies that
\begin{equation}\label{eq: contradiction}
B_{{i_0}}= \pm \beta \prod\limits_{i\not={i_0}}
d_{i,{\mathbf i}_0}^{b_{i,{\mathbf i}_0}}
= \pm \beta' \prod\limits_{i\not={i_0}}
d_{i,{\mathbf i}_1}^{b_{i,{\mathbf i}_1}},  \ \
 \  \
C_{{i_0}} = \pm\gamma\prod\limits_{i\not={i_0}} d_{i,{\mathbf i}_0}^{c_{i,{\mathbf i}_0}} = \pm\gamma'\prod\limits_{i\not={i_0}}
d_{i,{\mathbf i}_1}^{c_{i,{\mathbf i}_1}}\,.
\end{equation}
Since we assume that $B_{{i_0}},C_{{i_0}}\ne\const$,
\eqref{eq: contradiction} can occur only if
$\displaystyle B_{{i_0}}
=\pm\beta d_{s,{\mathbf i}_0}^{b_{s,{\mathbf i}_0}}$ and
$C_{{i_0}}=\pm\gamma d_{s,{\mathbf i}_0}^{c_{s,{\mathbf i}_0}}$.
I.e., $B_{{i_0}}$ and $C_{{i_0}}$ have a non-trivial common
factor, which contradicts our
assumptions. Thus, we have proved \eqref{eq:for A 1 general}.
\vskip5pt

\noindent By \eqref{claim 4}, the choice of $i_0$ implies that
$\sum a_{{i_0},{\mathbf i}}=\sum b_{{i_0},{\mathbf i}}
=\sum c_{{i_0},{\mathbf i}}>0$.
Since $b_{{i_0},{\mathbf i}_0}=0$, the latter implies that there
exists a multiindex ${\mathbf i}_2\ne {\mathbf i}_0$ such that
$b_{{i_0},{\mathbf i}_2}>0$. In the same way as we proved
\eqref{eq:for A 1 general}, we conclude that
\begin{equation}\label{eq:for B 1 general}
A_{{i_0}} = \pm \alpha' \prod\limits_{i\not={i_0}}
d_{i,{\mathbf i}_2}^{a_{i,{\mathbf i}_2}}\,,\qquad
B = B_{{i_0}} d_{{i_0},{\mathbf i}_2}^{b_{{i_0},{\mathbf i}_2}}
\,,\qquad
C_{{i_0}} = \pm\gamma'\prod\limits_{i\not={i_0}}
d_{i,{\mathbf i}_2}^{c_{i,{\mathbf i}_2}}\,.
\end{equation}
Comparing \eqref{eq:for A 1 general}
and \eqref{eq:for B 1 general}, we conclude that
$C_{{i_0}} = \pm\gamma d_{s,{\mathbf i}_0}^{c_{s,{\mathbf i}_0}}$
and
${\mathbf i}_0\cap{\mathbf i}_2\ne \varnothing$.
Pick $t_0\in {\mathbf i}_0\cap{\mathbf i}_2$, then $A_{t_0}=\const$.
Due to \eqref{claim 4},
$\deg A_{{i_0}} = \deg B_{{i_0}} =\deg C_{{i_0}}$; that is,
the polynomials $A_{t_0},B_{t_0},C_{t_0}$ are constant.
This completes the proof of $(*)$.
\vskip7pt

\noindent {\bf Step 2.} Let $i_0\in \{1,...,n\}$.
Let us prove the following statement:
\vskip5pt

$(**)$
\emph{Assume that the polynomials $A_{{i_0}},B_{{i_0}},C_{{i_0}}$ are
constant. Then there exist
a multiindex ${\mathbf i}=(i_1,\dots,i_{m-2})$, indices $j,s,l,t$,
and a complex number $\alpha\ne 0$ such that
all $i_0,i_1,\dots,i_{m-2},j,s,l,t$ are distinct and
$ A = \alpha d_{i_0,\mathbf{i},j,s} d_{i_0,\mathbf{i},l,t}$,
$ B = \alpha d_{i_0,\mathbf{i},j,l} d_{i_0,\mathbf{i},t,s}$ and
$ C=\alpha   d_{i_0,\mathbf{i},j,t} d_{i_0,\mathbf{i},s,l}$.}
\vskip5pt

{\it Proof.}
Denote
$$
A'= A_{{i_0}}\prod\delta_{{\mathbf i};1}^{a_{{i_0},{\mathbf i}}}\,,
\quad
B'= B_{{i_0}}\prod\delta_{{\mathbf i};1}^{b_{{i_0},{\mathbf i}}}\,,
\quad
C'= C_{{i_0}}\prod\delta_{{\mathbf i};1}^{c_{{i_0},{\mathbf i}}}\,.
$$
It is easily seen that $A'$, $B'$ and $C'$ are pairwise co-prime
on $\Bbb{C}^{(m-1)(n-1)}$ and do not vanish on the configuration space $\Cal{E}^{n-1}(\Bbb{C}^{m-1},gp)\subset \Bbb{C}^{(m-1)(n-1)}$. According to
\eqref{eq:needs Lin's theorem1 general}, $A'+B'+C'=0$.
The proof of $(**)$ is by induction on $m$.

Let $m=2$. Lemma $5.1$ of \cite{LinSphere} states that for
any three non-constant pairwise co-prime polynomials
$$
P=\mathrm{a}\prod_{i\ne j }(x_i-x_j)^{a_{i,j}}\,, \hfill
Q=\mathrm{b}\prod_{i\ne j }(x_i-x_j)^{b_{i,j}}\,, \hfill
R=\mathrm{c}\prod_{i\ne j }(x_i-x_j)^{c_{i,j}}
$$
in the variables $x_1,...,x_n$ which satisfy
the equation $A'+B'+C'=0$, there exist
distinct indices $j,s,l,t$ and
a complex number $\alpha \ne 0$ such that either
$$
P= \alpha (x_{j}-x_{s})\,, \
Q= \alpha (x_{s}-x_{l})\,, \
R= \alpha (x_{l}-x_{j})
$$
or
$$
P= \alpha (x_{j}-x_{s})(x_{l}-x_{t})\,, \
Q= \alpha (x_{j}-x_{l})(x_{t}-x_{s})\,, \
R= \alpha (x_{j}-x_{t})(x_{s}-x_{l})\,.
$$
This lemma applies to the
polynomials $A',B',C'$ in the variables $z_{1,2},...,z_{n,2}$.
It follows that there exist distinct indices $j,s,l,t$ and
a complex number $\alpha \ne 0$ such that either
$A'= \alpha (z_{j,2}-z_{s,2})$,
$B'= \alpha (z_{s,2}-z_{l,2})$,
$C'= \alpha (z_{l,2}-z_{j,2})$
or
$A'= \alpha (z_{j,2}-z_{s,2})(z_{l,2}-z_{t,2})$,
$B'= \alpha (z_{j,2}-z_{l,2})(z_{t,2}-z_{s,2})$,
$C'= \alpha (z_{j,2}-z_{t,2})(z_{s,2}-z_{l,2})$.
Thus, for our original polynomials $A,B,C$ we obtain that either
$$
A= \alpha d_{i_0,j,s}, \quad
B= \alpha d_{i_0,s,l}, \quad C= \alpha d_{i_0,l,j}
$$
or
$$
A= \alpha d_{i_0,j,s}d_{i_0,l,t}\,, \quad
B= \alpha d_{i_0,j,l}d_{i_0,t,s}\,, \quad
C= \alpha d_{i_0,j,t}d_{i_0,s,l}\,.
$$
It is easily seen that in the first case
$$
\aligned
A+B+C &=d_{i_0,j,s}+d_{i_0,s,l}+d_{i_0,l,j} \\
&=
z_{j,1}z_{s,2}-z_{s,1}z_{j,2}+z_{s,1}z_{l,2}-z_{l,1}
z_{s,2}+z_{l,1}z_{j,2}-z_{j,1}z_{l,2}\not\equiv 0\,.
\endaligned
$$
Thus, the equality $A+B+C=0$ can be fulfilled
only in the second case, which provides the base of induction.

Suppose that $(**)$ is fulfilled for some $m=k-1>1$
and let us prove that then the same is true for $m=k$.
Due to $(*)$ (from Step 1), the induction hypothesis
applies to the polynomials  $A',B',C'$; that is,
there exist a multiindex
${\mathbf i}=(i_1,...,i_{k-2})$, indices $j,s,l,t$,
and a complex number $\alpha\ne 0$ such that
all $i_1,...,i_{k-2},j,s,l,t$ are distinct and
$ A' = \alpha\cdot d_{\mathbf{i},j,s}\cdot d_{\mathbf{i},l,t}$,
$ B' = \alpha\cdot d_{\mathbf{i},j,l}\cdot d_{\mathbf{i},t,s}$ and
$ C' =\alpha\cdot d_{\mathbf{i},j,t}\cdot d_{\mathbf{i},s,l}$.
I.e., the original polynomials $A,B,C$ can be written as
$ A = \alpha\cdot d_{i_0,\mathbf{i},j,s}\cdot d_{i_0,\mathbf{i},l,t}$,
$ B = \alpha\cdot d_{i_0,\mathbf{i},j,l}\cdot d_{i_0,\mathbf{i},t,s}$
and
$ C =\alpha\cdot d_{i_0,\mathbf{i},j,t}\cdot d_{i_0,\mathbf{i},s,l}$.
This completes the justification of $(**)$
and proves Lemma \ref{Lm: likeABC for greater dimensions}.
\end{proof}

\subsection{Explicit description of all holomorphic
functions $\Cal{E}^n\to{\mathbb C}\setminus\{0,1\}$.}
\label{ss: Explicit description of all holomorphic functions }

\begin{Notation}\label{Not: set of all non-constant holomorphic functions Z to C-(0,1)}
For a complex space $Z$, we donote by
$L(Z)$ the set of all non-constant holomorphic functions
$Z\to{\mathbb C}\setminus\{0,1\}$.
\end{Notation}

\begin{Theorem}\label{Thm: holomorphic functions En(X) to C-(0,1)}
Let $X=\Bbb{C}^m$ or $X=\Bbb{CP}^m$. Then
$L(\Cal{E}^n(X,gp))=\dcr(\Cal{E}^n(X,gp))$\footnote{
See Definition \ref{Def: determinant cross ratio}.}.
\end{Theorem}

\begin{proof}
By Corollary \ref{Lm: determinant cross ratios omit 0 and 1},
it suffices to prove the inclusion
$$
L(\Cal{E}^n(X,gp))\subseteq \dcr(\Cal{E}^n(X,gp))\,.
$$
We follow \cite{Lin04b}. Let $\mu\in L(\Cal{E}^n(X,gp))$, that is,
$\mu\colon\Cal{E}^n(X,gp)\to{\mathbb C}\setminus\{0,1\}$
is a holomorphic function.

First, let $X=\Bbb{C}^m$.
It follows from Big Picard Theorem that $\mu$
is a regular function on $\Cal{E}^n(\Bbb{C}^m,gp)$;
hence it is a rational function
on $(\Bbb{C}^m)^n$ and there are co-prime polynomials
$A,B\in\Bbb{C}[(\Bbb{C}^m)^n]$
that do not vanish on $\Cal{E}^n(\Bbb{C}^m,gp)$ and such that
$\mu=-A/B$. The function $1-\mu=(A+B)/B$ also
omits the values $0$, $1$. The polynomials $A$, $B$ and
$C=-B-A$ are pairwise co-prime, do not vanish on
$\Cal{E}^n(\Bbb{C}^m,gp)$ and satisfy $A+B+C=0$.
Lemma \ref{Lm: likeABC for greater dimensions} applies to the
latter three polynomials and shows that
$\mu=-A/B=-d_{\mathbf{i},j,k} d_{\mathbf{i},l,s}/
d_{\mathbf{i},j,l} d_{\mathbf{i},s,k}=
e_{{\mathbf i};j,k,l,s}$ for appropriate ${\mathbf i},j,k,l,s$.

When $X=\Bbb{CP}^m$, we restrict $\mu$ from
$\Cal{E}^n(\Bbb{CP}^m,gp)$
to $\Cal{E}^n(\Bbb{C}^m,gp)$ and apply the above result,
which leads to the desired conclusion.
\end{proof}

\noindent The following lemma is a known fact of the
classical invariant theory (for small dimension it was discovered
by A. F. M{\" o}bius, \cite{Mobius27}, especially, Part 2).
The proof may be extracted from
\cite{Weyl46}, Section 2.14 (especially, Theorem 2.14.A).
However, the exposition in \cite{Weyl46} is rather
complicated; for the reader's convenience,
we give here an independent proof.

\begin{Lemma}\label{Lm: Sum of determinant is 0}
Let $\mathbf{i}=(i_1,\dots,i_{m-1})$ be a multiindex,
and $j,k,l,s$ be indices such that
all $i_1,\dots,i_{m-1},j,k,l,s$ are distinct. Then
$d_{\mathbf{i},j,k} d_{\mathbf{i},l,s}+
 d_{\mathbf{i},j,l} d_{\mathbf{i},s,k}+
 d_{\mathbf{i},j,s} d_{\mathbf{i},k,l}=0$.
\end{Lemma}

\begin{proof}
By Lemma \ref{Lm: determinant cross ratios omit 0 and 1},
the function $1-e_{{{\mathbf i}};j,k,l,s}$
omits the values $0$ and $1$ on $\Cal{E}^n$;
thus, by Theorem \ref{Thm: holomorphic functions En(X) to C-(0,1)},
it is a determinant cross ratio, say
$e_{{{\mathbf i}'};j',k',l',s'}$.
That is,
$$
1-e_{{{\mathbf i}};j,k,l,s}=
1- \frac{d_{{\mathbf i},j,k}d_{{\mathbf i},l,s}}
{d_{{\mathbf i},j,l}d_{{\mathbf i},k,s}} =
\frac{d_{{\mathbf i},j,l}d_{{\mathbf i},k,s}-
d_{{\mathbf i},j,k}d_{{\mathbf i},l,s}}
{d_{{\mathbf i},j,l}d_{{\mathbf i},k,s}} = e_{{{\mathbf i}'};j',k',l',s'}=\frac{d_{{\mathbf i}',j',k'}d_{{\mathbf i}',l',s'}}
{d_{{\mathbf i}',j',l'}d_{{\mathbf i}',k',s'}}\,.
$$
Since determinant polynomials are irreducible (see Lemma
\ref{cl: irreducibility}) and the polynomials
$d_{{\mathbf i}',j',k'}d_{{\mathbf i}',l',s'}$
and $d_{{\mathbf i}',j',l'}d_{{\mathbf i}',k',s'}$ are co-prime, it follows that there is a complex number $c\ne 0$ such that
\begin{equation}\label{eqality between determinants}
d_{{\mathbf i},j,l}d_{{\mathbf i},k,s}-
d_{{\mathbf i},j,k}d_{{\mathbf i},l,s}
=c d_{{\mathbf i}',j',k'}d_{{\mathbf i}',l',s'}
\end{equation}
and
\begin{equation}\label{eqality between determinants1}
d_{{\mathbf i},j,l}d_{{\mathbf i},k,s}
=c d_{{\mathbf i}',j',l'}d_{{\mathbf i}',k',s'}\,.
\end{equation}
By definition of a determinant cross ratio
and Lemma \ref{cl: irreducibility}, the polynomials
$d_{{\mathbf i},j,l}d_{{\mathbf i},k,s}$ and
$d_{{\mathbf i},j,k}d_{{\mathbf i},l,s}$ are co-prime.
From \eqref{eqality between determinants} it follows that
$d_{{\mathbf i},j,l}d_{{\mathbf i},k,s}$,
$d_{{\mathbf i},j,k}d_{{\mathbf i},l,s}$ and
$d_{{\mathbf i}',j',k'}d_{{\mathbf i}',l',s'}$ are pairwise
co-prime.
By Lemma \ref{Lm: likeABC for greater dimensions},
equality \eqref{eqality between determinants} implies
that
\begin{equation}\label{eqality between determinants2}
c d_{{\mathbf i}',j',k'}d_{{\mathbf i}',l',s'}=
d_{{\mathbf i},j,s}d_{{\mathbf i},k,l}\,.
\end{equation}
By \eqref{eqality between determinants},
$d_{{\mathbf i},j,l}d_{{\mathbf i},k,s}-
d_{{\mathbf i},j,k}d_{{\mathbf i},l,s}=
d_{{\mathbf i},j,s}d_{{\mathbf i},k,l}$;
that is,
$$
d_{\mathbf{i},j,k} d_{\mathbf{i},l,s}+
 d_{\mathbf{i},j,l} d_{\mathbf{i},s,k}+
 d_{\mathbf{i},j,s} d_{\mathbf{i},k,l}=0\,.
$$
%
%
\end{proof}

\begin{Corollary}\label{sum of determinant cross ratios}
$e_{{\mathbf i};j,k,l,s}+e_{{\mathbf i};j,s,l,k}=1$.
\end{Corollary}

\begin{proof}
\noindent Indeed, by the above lemma,
\begin{equation*}
\aligned
e_{{\mathbf i};j,k,l,s}(q)+e_{{\mathbf i};j,s,l,k}(q)-1&=
\frac{d_{{\mathbf i},j,k}(q)}{d_{{\mathbf i},j,l}(q)}:
\frac{d_{{\mathbf i},k,s}(q)}{d_{{\mathbf i},l,s}(q)}
+
\frac{d_{{\mathbf i},j,s}(q)}{d_{{\mathbf i},j,l}(q)}:
\frac{d_{{\mathbf i},s,k}(q)}{d_{{\mathbf i},l,k}(q)}-1\\
&=\frac{d_{{\mathbf i},j,k}(q)}{d_{{\mathbf i},j,l}(q)}:
\frac{d_{{\mathbf i},k,s}(q)}{d_{{\mathbf i},l,s}(q)}
+
\frac{d_{{\mathbf i},j,s}(q)}{d_{{\mathbf i},j,l}(q)}:
\frac{d_{{\mathbf i},k,s}(q)}{d_{{\mathbf i},k,l}(q)}-1
\\
&=\frac{
d_{{\mathbf i},j,k}(q)d_{{\mathbf i},l,s}(q)+
d_{{\mathbf i},j,s}(q)d_{{\mathbf i},k,l}(q)+
d_{{\mathbf i},j,l}(q)d_{{\mathbf i},s,k}(q)
}
{d_{{\mathbf i},j,l}(q)d_{{\mathbf i},k,s}(q)}\\
& = 0\,.
\endaligned
\end{equation*}
\end{proof}

\noindent In the following lemma we establish some
simple relations between determinant cross ratios.

\begin{Lemma}
\label{Lm: 3 facts about strong equivariance multi dim}
$a)$ $e_{{\mathbf i};t,k,r,s} = e_{{\mathbf i};j,k,r,s}/
e_{{\mathbf i};j,k,r,t}$;\\
$b)$ $\displaystyle e_{{\mathbf i};j,k,t,s} =
1-\frac{1-e_{{\mathbf i};j,k,r,s} }{1- e_{{\mathbf i};j,k,r,t}}$;\\
$c)$ $\displaystyle e_{{\mathbf i};j,t,r,s} =
\left(1-
\frac{1-(e_{{\mathbf i};j,k,r,s})^{-1}}
{1-(e_{{\mathbf i};j,k,r,t})^{-1}}\right)^{-1}$;\\
$d)$ $
e_{{\mathbf j},i;j,k,r,s}=e_{{\mathbf j},j;i,k,r,s}
e_{{\mathbf j},s;j,k,r,i}=
e_{{\mathbf j},k;j,i,r,s}e_{{\mathbf j},r;j,k,i,s}
$.
\end{Lemma}

\begin{proof}
$(a)$ follows directly form the definition of
determinant cross ratios.

$b)$ By Corollary \ref{sum of determinant cross ratios},
$$
\aligned
1-e_{{\mathbf i};j,k,r,s}=e_{{\mathbf i};j,s,r,k} =
e_{{\mathbf i};r,k,j,s} \,,\\
1-e_{{\mathbf i};j,k,r,t}=e_{{\mathbf i};r,k,j,t} =
e_{{\mathbf i};r,k,j,t}\,.
\endaligned
$$
In view of $(a)$ and
Corollary \ref{sum of determinant cross ratios}, this implies
$$
\frac{1-e_{{\mathbf i};j,k,r,s} }{1- e_{{\mathbf i};j,k,r,t}}=
\frac{e_{{\mathbf i};r,k,j,s} }{e_{{\mathbf i};r,k,j,t}}=
e_{{\mathbf i};t,k,j,s}=1-e_{{\mathbf i};j,k,t,s}\,.
$$

$c)$ Clearly,
$(e_{{\mathbf i};j,k,r,s})^{-1}=e_{{\mathbf i};j,r,k,s}$
and $(e_{{\mathbf i};j,k,r,t})^{-1}=e_{{\mathbf i};j,r,k,t}$.
By $(b)$,
$$
\left(1-
\frac{1-(e_{{\mathbf i};j,k,r,s})^{-1}}
{1-(e_{{\mathbf i};j,k,r,t})^{-1}}\right)=
1-\frac{1-e_{{\mathbf i};j,r,k,s} }{1- e_{{\mathbf i};j,r,k,t}}
= e_{{\mathbf i};j,r,t,s}\,,
$$
since $(e_{{\mathbf i};j,r,t,s})^{-1}= e_{{\mathbf i};j,t,r,s}$,
this implies $(c)$.

$d)$ From the definition of a determinant cross ratio
we have
\begin{equation}\label{eq: some long equation}
\aligned
e_{{\mathbf j},i;j,k,r,s} / e_{{\mathbf j},j;i,k,r,s} &=
\frac{d_{{\mathbf j},i,j,k}}{d_{{\mathbf j},i,j,r}}
  \frac{d_{{\mathbf j},i,r,s}}{d_{{\mathbf j},i,k,s}}
  \frac{d_{{\mathbf j},j,i,r}}{d_{{\mathbf j},j,i,k}}
  \frac{d_{{\mathbf j},j,k,s}}{d_{{\mathbf j},j,r,s}}\\
  &=\frac{d_{{\mathbf j},i,r,s}}{d_{{\mathbf j},i,k,s}}
\frac{d_{{\mathbf j},j,k,s}}{d_{{\mathbf j},j,r,s}}
= \frac{d_{{\mathbf j},s,j,k}}{d_{{\mathbf j},s,j,r}}
\frac{d_{{\mathbf j},s,i,r}}{d_{{\mathbf j},s,i,k}}
= e_{{\mathbf j},s;j,k,r,i}\,.
\endaligned
\end{equation}
Clearly, $e_{{\mathbf j},i;j,k,r,s}= e_{{\mathbf j},i;k,j,s,r}$
and
$e_{{\mathbf j},k;j,i,r,s} = e_{{\mathbf j},k;i,j,s,r}$.
By \eqref{eq: some long equation},
$$
e_{{\mathbf j},i;j,k,r,s}/e_{{\mathbf j},k;j,i,r,s} =
e_{{\mathbf j},i;k,j,s,r}/e_{{\mathbf j},k;i,j,s,r} =
e_{{\mathbf j},r;k,j,s,i} = e_{{\mathbf j},r;j,k,i,s}\,.
$$
\end{proof}

\section{Simplicial complex of holomorphic functions
$\Cal{E}^n\to{\mathbb C}\setminus\{0,1\}$}
\label{sec: Symplicial complex of holomorphic functions En to C-0 1}

\subsection{Simplicial complex of holomorphic functions omitting two values}
\label{ss: Symplicial complex of holomorphic functions omitting two values}
It was shown in \cite{LinSphere} that
the set $L(Z)$ of all non-constant holomorphic functions
$Z\to{\mathbb C}\setminus\{0,1\}$ on a complex space $Z$
may be endowed with a natural structure of a simplicial complex
$L_{\vartriangle}(Z)$ and the correspondence $Z\mapsto L_{\vartriangle}(Z)$
has some properties of a contravariant functor from the category
of complex spaces and holomorphic mappings to the category of
simplicial complexes and simplicial mappings.
In this section, we apply this Lin's construction for $Z=\Cal{E}^n$
and study some properties of the complex $L_{\vartriangle}(\Cal{E}^n)$.
First, we recall the definition of the complex $L_{\vartriangle}(Z)$.

\begin{Definition}\label{Def: divisibility and simplices}
Let $Z$ be a complex space and $L(Z)$ be the set of all
non-constant holomorphic functions
$Z\to{\mathbb C}\setminus\{0,1\}$.

For $\mu,\nu\in L(Z)$, we say that $\nu$
is {\em a proper divisor of} $\mu$ and write $\nu\mid\mu$
if the quotient $\lambda=\mu:\nu\in L(Z)$, i. e.,
$\lambda\neq\const$ and $\lambda(z)\ne 1$ for
all $z\in Z$; otherwise, we write $\nu\nmid\mu$.
Clearly, $\nu\mid\mu$ is equivalent to $\mu\mid\nu$.

A non-empty ordered subset $\Delta=\{\mu_0,...,\mu_s\}
\subseteq L(Z)$ is said to be a simplex of dimension
$s$ with vertices $\mu_0,...,\mu_s$ if $\mu_i\mid\mu_j$
for each pair of distinct $i,j$.
Evidently, a non-empty subset of a simplex is a simplex,
that is, we obtain a well-defined simplicial complex $L_{\vartriangle}(Z)$
with the set of vertices $L(Z)$.\footnote
{It is well-known that if $Z$ is a quasi-projective algebraic variety,
then the set $L(Z)$ consists of a finite number of regular
functions and the complex $L_{\vartriangle}(Z)$
is finite (it can be empty); however, we do not need the latter
statement since we already know the explicit description of
$L(\Cal{E}^n)$.}

If $Z$ is a quasi-projective algebraic variety,
then the set $L(Z)$ consists of a finite number of regular
functions and the complex $L_{\vartriangle}(Z)$
is finite (it can be empty).

A holomorphic map $f\colon Z\to Y$ of complex spaces
induces the homomorphism
$f^*\colon{\mathcal O}(Y)\to{\mathcal O}(Z)$
of the algebras of holomorphic functions defined, as usual, by
$f^*(\lambda)=\lambda\circ f$, $\lambda\in{\mathcal O}(Y)$.
Let $\lambda\in L(Y)$; suppose that $f^*(\lambda)\ne\const$
(this is certainly the case whenever $f$ is a dominant map
of irreducible quasi-projective varieties).
Then the map of the vertices
$f^*\colon L(Y)\ni\lambda\mapsto\lambda\circ f\in L(Z)$
induces the simplicial map $f^*\colon L_{\vartriangle}(Y)
\to L_{\vartriangle}(Z)$ whose restriction to each
simplex $\Delta\subseteq L(Y)$ is injective and
preserves dimensions of simplices.
\hfill $\bigcirc$
\end{Definition}

\begin{Remark}\label{Rmk: degradation}
Let ${\mathbf i}=(i_1,\dots,i_{m-1})$ and
$$
\mu=e_{{\mathbf i};j,k,s,t}=
\frac{d_{{\mathbf i},j,k}}{d_{{\mathbf i},j,s}}:
\frac{d_{{\mathbf i},k,t}}{d_{{\mathbf i},s,t}}\,.
$$
Let ${\mathbf e}_1,...,{\mathbf e}_m$ be the standard basis in
$\Bbb{C}^m$, ${\mathbf u}_1=0$ and
${\mathbf u}_p=\sum\limits_{j=m-p+2}^{m}
{\mathbf e}_j$ for $p=2,...,m-1$.
Set
$$
L = \{(v_1,...,v_n)\in (\Bbb{C}^m)^n\,| v_{i_p} ={\mathbf u}_p \
\text{\rm for} \ p=1,...,m-1\}\,.
$$
For $q=\left(\begin{array}{c}q_1\\ \vdots\\ q_n\end{array}
\right)\in (\Bbb{C}^m)^n \cap L$
with all $q_r=(z_{r,1},...,z_{r,m})\in{\mathbb C}^m$,
the restriction of any determinant polynomial
$d_{{\mathbf i},l,r}$ to the subspace $L$ is reduced
to a certain determinant of order 2 and may be computed
as $(-1)^m (z_{l,1}z_{r,2}-z_{r,1}z_{l,2})$.
Consequently,
the restriction $\mu|_L$ of the function $\mu$
to the subspace $L$ may be written as
\begin{equation}\label{eq: usual cross ratio}
\mu|_L(q)=
\frac{(z_{j,1}/z_{j,2})-(z_{k,1}/z_{k,2})}
{(z_{j,1}/z_{j,2})-(z_{s,1}/z_{s,2})}
:
\frac{(z_{k,1}/z_{k,2})- (z_{t,1}/z_{t,2})}
{(z_{s,1}/z_{s,2})-(z_{t,1}/z_{t,2})}\,.
\end{equation}
Thus, $\mu|_L$ is the cross ratio of the four quantities
$a=z_{j,1}/z_{j,2}$, $b=z_{k,1}/z_{k,2}$,
$c=z_{j,1}/z_{j,2}$ and $d=z_{t,1}/z_{t,2}$,
which may be treated as four distinct
points in ${\mathbb{CP}}^1$ whenever $q\in{\mathcal E}^n\cap L$.
If $\mu$ and $\mu'$ are determinant cross ratios,
$\mu|\mu'$ and $\mathbf{i} = \supp_{\text{ess}}\mu
=\supp_{\text{ess}}\mu'$, then
$\supp_{\text{ess}}(\mu:\mu')
= \supp_{\text{ess}}\mu=\supp_{\text{ess}}\mu'$.
In view of the above,
the restrictions of $\mu$, $\mu'$ and $\mu:\mu'$
to $L$ are usual cross ratios and, moreover,
$(\mu|_L):(\mu'|_L)=(\mu:\mu')|_L$.
\hfill $\bigcirc$
\end{Remark}

\noindent
In view of the previous remark, certain results about  usual cross
ratios apply to the determinant cross ratios, as well.
In particular, we will use the following simple lemma
(compare to Lemma 5.7, \cite{LinSphere}):

\begin{Lemma}\label{Lm: simplices of cross ratios}
If the ratio $\mu:\mu'$ of two
cross ratios
$$
\mu=\frac{q_j-q_k}{q_j-q_s}:\frac{q_k-q_t}{q_s-q_t}\ \
\text{\rm and} \ \
\mu'=\frac{q_{j'}-q_{k'}}{q_{j'}-q_{s'}}:
\frac{q_{k'}-q_{t'}}{q_{s'}-q_{t'}}
$$
is a cross ratio of certain four of the eight variables
$q_j,q_k,q_s,q_t,q_{j'},q_{k'},q_{s'},q_{t'}$
then $\card(\{j,k,s,t\}\cap\{j',k',s',t'\})=3$ and
$\mu'$ is obtained from $\mu$ by replacing of one
of the variables $q_j,q_k,q_s,q_t$ with some $q_m$,
where $m\ne j,k,s,t$.
\hfill $\bigcirc$
\end{Lemma}

\noindent We need also the following technical lemma.

\begin{Lemma}\label{Lm: m-generalized cross ratios lemma part a}
Let $X$ be either $\Bbb{C}^m$ or $\Bbb{CP}^m$, $n>m+2$, and
let $\mu=e_{\mathbf{i};j,k,r,s}$ be a proper divisor of a
determinant cross ratio $\mu'$. Then
$\card(\supp\limits_{\text{ess}}\mu \cap
\supp\limits_{\text{ess}}\mu')\ge m-2$.
If $\supp\mu'\not=\supp\mu$, then
$\card(\supp\mu\cap\supp\mu')=m+2$,
$\supp_{\text{ess}}\mu =\supp_{\text{ess}}\mu'$, and
$\mu'$ is one of the functions
$e_{{\mathbf i};j,k,r,t}$, $e_{{\mathbf i};k,j,s,t}$,
$e_{{\mathbf i};r,s,j,t}$, $e_{{\mathbf i};s,r,k,t}$.
\end{Lemma}

\begin{proof}
Set
\begin{equation*}
 \mu = e_{{{\mathbf i}};j,k,r,s}=
    \frac{d_{{\mathbf i},j,k}}{d_{{\mathbf i},j,r}}:
  \frac{d_{{\mathbf i},k,s}}{d_{{\mathbf i},r,s}}\,,
\quad \mu' = e_{{\mathbf i}';j',k',r',s'}=
    \frac{d_{{\mathbf i}',j',k'}}{d_{{\mathbf i}',j',r'}}:
  \frac{d_{{\mathbf i}',k',s'}}{d_{{\mathbf i}',r',s'}}\,.
 \end{equation*}
Let ${\mathbf j}={\mathbf i}\cap{\mathbf i}'$.
Since $\mu|\mu'$, from Lemma \ref{cl: irreducibility}
it follows that
\begin{equation}\label{eq: about support}
{\mathbf i}, {\mathbf i}'\in (\supp\mu\cap\supp\mu')
\subset \supp\mu\,.
\end{equation}
Let us show that $\card{\mathbf j}\ne m-5,m-4,m-3$
and hence $\card{\mathbf j}\ge m-2$.

If $\card{\mathbf j} = m-5$, then
by \eqref{eq: about support},
${\mathbf i} = ({\mathbf j},j',k',r',s')$
and ${\mathbf i}' = ({\mathbf j},j,k,r,s)$.
Thus,
\begin{gather*}
 \mu = e_{{\mathbf j},j',k',r',s';j,k,r,s}=
  \frac{d_{{\mathbf j},j' ,k',r',s',j,k}}
 {d_{{\mathbf j},j',k',r',s',j,r}}:
\frac{d_{{\mathbf j},j',k',r',s',k,s}}
{d_{{\mathbf j},j',k',r',s',r,s}}\,,\\
\mu' = e_{{\mathbf j},j,k,r,s;j',k',r',s'}=
\frac{d_{{\mathbf j},j,k,r,s,j',k'}}
{d_{{\mathbf j},j,k,r,s,j',r'}}:
\frac{d_{{\mathbf j},j,k,r,s,k',s'}}
{d_{{\mathbf j},j,k,r,s,r',s'}}\,.
 \end{gather*}
Since $\mu/\mu'$ is a determinant cross ratio,
and determinant polynomials are irreducible
(Lemma \ref{cl: irreducibility}), $\card(\{j',k',r',s'\}
\cap \{j,k,r,s\}) > 0$, and we come
to $\card({\mathbf i}\cap{\mathbf i}')> m-5$,
which contradicts our assumption.

For $\card({\mathbf i}\cap{\mathbf i}')=m-4$,
in view of \eqref{eq: about support},
we can assume that
${\mathbf i} = ({\mathbf j},j',k',r')$ and
${\mathbf i}' = ({\mathbf j},j,k,r)$.
Thus,
\begin{gather*}
 \mu = e_{{\mathbf j},j',k',r';j,k,r,s}=
 \frac{d_{{\mathbf j},j',k',r',j,k}
 }{d_{{\mathbf j},j',k',r',j,r}}:
\frac{d_{{\mathbf j},j',k',r',k,s}
}{d_{{\mathbf j},j',k',r',r,s}}\,,\\
\mu' = e_{{\mathbf j},j,k,r;j',k',r',s'}=
\frac{d_{{\mathbf j},j,k,r,j',k'}}
{d_{{\mathbf j},j,k,r,j',r'}}:
\frac{d_{{\mathbf j},j,k,r,k',s'}}
{d_{{\mathbf j},j,k,r,r',s'}}\,.
 \end{gather*}
Since $\mu/\mu'$ is a determinant cross ratio
and determinant polynomials are irreducible
(Lemma \ref{cl: irreducibility}),
$\card(\{j',k',r'\} \cap \{j,k,r\})> 0$,
a contradiction.

For $\card({\mathbf i}\cap{\mathbf i}')=m-3$,
in view of \eqref{eq: about support},
we can assume that either
\begin{gather*}
 \mu = e_{{\mathbf j},j',k';j,k,r,s}=
 \frac{d_{{\mathbf j},j',k',j,k}
 }{d_{{\mathbf j},j',k',j,r}}:
\frac{d_{{\mathbf j},j',k',k,s}
}{d_{{\mathbf j},j',k',r,s}}\,,\\
\mu' = e_{{\mathbf j},j,k;j',k',r',s'}=
\frac{d_{{\mathbf j},j,k,j',k'}}
{d_{{\mathbf j},j,k,j',r'}}:
\frac{d_{{\mathbf j},j,k,k',s'}}
{d_{{\mathbf j},j,k,r',s'}}
 \end{gather*}
or
\begin{gather*}
 \mu = e_{{\mathbf j},j',k';j,k,r,s}=
 \frac{d_{{\mathbf j},j',k',j,k}
 }{d_{{\mathbf j},j',k',j,r}}:
\frac{d_{{\mathbf j},j',k',k,s}
}{d_{{\mathbf j},j',k',r,s}}\,,\\
\mu' = e_{{\mathbf j},j,r;j',k',r',s'}=
\frac{d_{{\mathbf j},j,r,j',k'}}
{d_{{\mathbf j},j,r,j',r'}}:
\frac{d_{{\mathbf j},j,r,k',s'}}
{d_{{\mathbf j},j,r,r',s'}}
 \end{gather*}
or
\begin{gather*}
 \mu = e_{{\mathbf j},j',k';j,k,r,s}=
 \frac{d_{{\mathbf j},j',k',j,k}
 }{d_{{\mathbf j},j',k',j,r}}:
\frac{d_{{\mathbf j},j',k',k,s}
}{d_{{\mathbf j},j',k',r,s}}\,,\\
\mu' = e_{{\mathbf j},j,s;j',k',r',s'}=
\frac{d_{{\mathbf j},j,s,j',k'}}
{d_{{\mathbf j},j,s,j',r'}}:
\frac{d_{{\mathbf j},j,s,k',s'}}
{d_{{\mathbf j},j,s,r',s'}}
 \end{gather*}
or
 \begin{gather*}
 \mu = e_{{\mathbf j},j',s';j,k,r,s}=
 \frac{d_{{\mathbf j},j',s',j,k}
 }{d_{{\mathbf j},j',s',j,r}}:
\frac{d_{{\mathbf j},j',s',k,s}
}{d_{{\mathbf j},j',s',r,s}}\,,\\
\mu' = e_{{\mathbf j},j,s;j',k',r',s'}=
\frac{d_{{\mathbf j},j,s,j',k'}}
{d_{{\mathbf j},j,s,j',r'}}:
\frac{d_{{\mathbf j},j,s,k',s'}}
{d_{{\mathbf j},j,s,r',s'}}\,.
 \end{gather*}

\noindent Since $\mu/\mu'$ is a determinant cross ratio
and determinant polynomials are irreducible
(Lemma \ref{cl: irreducibility}), these cases are impossible,
a contradiction.
This completes the proof of the first part
of the lemma.
\vskip5pt

\noindent Suppose now that  $\supp\mu \ne \supp{\mu'}$.
Since we have already proved that
$m>\card{\mathbf j}\ge m-2$,
we need to show that
$\card{\mathbf j}\ne m-2$.

Suppose to the contrary that
$\card{\mathbf j}= m-2$. Then,
without loss of generality, we may assume that
$\mu = e_{{\mathbf j},j';j,k,r,s}$ and
$\mu' =e_{{\mathbf j},j;j',k',r',s'}$.
Therefore,
\begin{equation*}
\aligned
 \mu : \mu'=&\left(\frac{d_{{\mathbf j},j',j,k}}
 {d_{{\mathbf j},j',j,r}}:
  \frac{d_{{\mathbf j},j',k,s}}
  {d_{{\mathbf j},j',r,s}}\right):
 \left(\frac{d_{{\mathbf j},j,j',k'}}
 {d_{{\mathbf j},j,j',r'}}:
 \frac{d_{{\mathbf j},j,k',s'}}
{d_{{\mathbf j},j,r',s'}}\right)\\
=&\frac{d_{{\mathbf j},j',j,k}
d_{{\mathbf j},j',r,s}
d_{{\mathbf j},j,j',r'}
d_{{\mathbf j},j,k',s'}}
{d_{{\mathbf j},j,j',k'}
d_{{\mathbf j},j',k,s}
d_{{\mathbf j},j',j,r}
d_{{\mathbf j},j,r',s'}}
\,.
\endaligned
\end{equation*}
Since all determinant polynomials are irreducible
(Lemma \ref{cl: irreducibility}), the latter quotient
is a determinant cross ratio if and only if
$k' = k$, $r' = r$ and $s'= s$,
which may happen if and only if $\supp\mu = \supp{\mu'}$,
a contradiction.

We are left with the case
$\card {\mathbf j}=m-1$, that is,
${\mathbf i}={\mathbf i}'$.
Let $L$ be as in Remark \ref{Rmk: degradation}; then
the restrictions $\mu|_L$ and $\mu'|_L$ of $\mu$ and $\mu'$
to $L$ are usual cross ratios of the variables
$p_j,p_k,p_r,p_s$ and $p_{j'},p_{k'},p_{r'},p_{s'}$,
respectively, and $(\mu|_L):(\mu'|_L)$ is such
a cross ratio, as well. By Lemma \ref{Lm: simplices of cross ratios},
$\card(\{j,k,r,s\} \cap \{j',k',r',s'\}) =3$
and the ordered set $\{j',k',r',s'\}$ is obtained
by replacing of one
of the indices in the ordered set $\{j,k,r,s\}$ with some index $t$,
where $t\ne j,k,r,s$ (up to a Kleinian permutation).
Consequently,
$\card(\supp\mu \cap \supp\mu')
=\card {{\mathbf i}} + 3 = m+2$ and
$\mu'$ is one of the functions
$e_{{\mathbf i};j,k,r,t}$, $e_{{\mathbf i};k,j,s,t}$,
$e_{{\mathbf i};r,s,j,t}$, $e_{{\mathbf i};s,r,k,t}$.
\end{proof}

\subsection{${\mathbf S}(n)$ action in
$L_{\vartriangle}(\mathcal{E}^n)$}
\label{S(n) action on L(En)}
The ${\mathbf S}(n)$ action in ${\mathcal E}^n$
induces an ${\mathbf S}(n)$ action on the set
$L({\mathcal E}^n)$ of all non-constant holomorphic functions
${\mathcal E}^n\to{\mathbb C}\setminus\{0,1\}$.
Of course, this is the action on the set of all
determinant cross ratios that we dealt with in
Lemma \ref{Rmk: transitivity of cross ratios}; thus,
it is transitive.

If $\mu,\nu\in L({\mathcal E}^n)$
and $\sigma\in{\mathbf S}(n)$, then the relations
$\mu|\nu$ and $(\sigma\mu)|(\sigma\nu)$ are equivalent.
Therefore, the above action induces a simplicial
dimension preserving ${\mathbf S}(n)$ action on the
complex $L_{\vartriangle}(\mathcal{E}^n)$; that is,
${\mathbf S}(n)$ acts on the set of all simplices of any
fixed dimension $s$. Our nearest goal is to describe
the orbits of the latter action;
in particular, we shall prove that on the set of all simplices
of any positive dimension it has exactly two orbits.
\vskip5pt

\noindent According to Notation \ref{Nt: first m members},
${\mathbf m}(\hat s) = (1,...,\hat{s}, ...,m)$;
sometimes we write $\widehat{\mathbf m}$
instead of ${\mathbf m}(\hat m)$.

\begin{Definition}\label{def: normal simplices}
The $t$-simplex
$$
\nabla_1^t=
\{e_{\widehat{\mathbf m};m,m+1,m+2,m+3};
e_{\widehat{\mathbf m};m,m+1,m+2,m+4};
...;e_{\widehat{\mathbf m};m,m+1,m+2,m+3+t}\}
$$
is called the {\em normal $t$-simplex of the first type};
such simplices do exist for $0\le t \le n-m-3$.
We say that a $t$-simplex is of the {\em first type}
if it belongs to the ${\mathbf S}(n)$ orbit of the
normal simplex $\nabla_1^t$.

The $t$-simplex
$$
\aligned
\nabla_2^t=\{e_{{\mathbf m}(\widehat{m-t});m-t,m+1,m+2,m+3};
&e_{{\mathbf m}(\widehat{m-t+1});m-t+1,m+1,m+2,m+3};\\
&...;
e_{{\mathbf m}(\hat m);m,m+1,m+2,m+3}\}
\endaligned
$$
is called the {\em normal $t$-simplex of the second type};
such simplices do exist for $0\le t \le m-1$.
We say that a $t$-simplex is of the {\em second type}
if it belongs to the ${\mathbf S}(n)$ orbit of the
normal simplex $\nabla_2^t$.

Notice that for $t>0$ the simplices
$\nabla_1^t$ and $\nabla_2^t$ belong to
different orbits of the ${\mathbf S}(n)$ action
on the set of all $t$-simpices.
\hfill $\bigcirc$
\end{Definition}

\begin{Remark}
\label{rm: simplices involving the same set of vector variables}
If $\Delta_1=\{\mu_0,...,\mu_t\}$ and
$\Delta_2=\{\nu_0,...,\nu_t\}$ are simplices of the same type
and the corresponding sets of functions $\mu$'s and $\nu$'s
both involve only the vector variables
$q_{i_1},...,q_{i_r}$, then $\Delta_1$ may be carried to
$\Delta_2$ by a permutation
$\sigma\in{\mathbf S}(\{i_1,...,i_r\})\subset{\mathbf S}(n)$.
\hfill $\bigcirc$
\end{Remark}


\noindent The following lemma shows that the ${\mathbf S}(n)$
action on the set of all simplices of dimension $t>0$
has at most two orbits and each of these orbits
consists of all simplices of the same type.
Notice that, by Theorem \ref{Thm: holomorphic functions En(X) to C-(0,1)},
any vertex of the complex
$L_{\vartriangle}(\mathcal{E}^n)$ is a determinant cross ratio.

\begin{Lemma}\label{Lm: m-generalized cross ratios lemma}
Let $X$ be either $\Bbb{C}^m$ or $\Bbb{CP}^m$ and $n>m+2$.

$a)$ Let
$\Delta=\{\mu_s\}_{s=1}^{l+1}\in L_{\vartriangle}({\mathcal E}^n)$
be an $l$-dimensional simplex and let
$\supp\mu_1=\supp\mu_2=...=\supp\mu_{l+1}$. Then
$\card(\supp\limits_{\text{ess}}\mu_s
\cap\supp\limits_{\text{ess}}\mu_t)=m-2$
for all $t\neq s$ and $l\le m-1$. Moreover, the simplex $\Delta$
is of the second type.

$b)$
Let
$\Delta=\{\mu_s\}_{s=1}^{l+1}\in
L_{\vartriangle}({\mathcal E}^n)$
be an $l$-dimensional simplex. If $\supp \mu_{s_0}
\ne\supp \mu_{t_0}$ for some $s_0 \ne t_0$,
then $\supp \mu_s \ne\supp \mu_t$ and
$\supp\limits_{\text{ess}} \mu_s=\supp\limits_{\text{ess}}\mu_t$
for all $s\ne t$. Moreover, the simplex $\Delta$
is of the first type.

$c)$ $\dim  L_{\vartriangle}({\mathcal E}^n(X))
= \max \{n-(m+3),m-1\}$.
\end{Lemma}

\begin{proof}
$a)$
Suppose, to the contrary, that
$\card(\supp\limits_{\text{ess}}\mu_s \cap
\supp\limits_{\text{ess}}\mu_t)\ne m-2$ for certain $s\ne t$.
By Lemma \ref{Lm: m-generalized cross ratios lemma part a},
this means that
$\card(\supp_{\text{ess}}\mu_s \cap
\supp_{\text{ess}}\mu_t)> m-2$.
Since $\card \supp_{\text{ess}}\mu_s=
\card \supp_{\text{ess}}\mu_t=m-1$, it follows that
$\supp_{\text{ess}}\mu_s=\supp_{\text{ess}}\mu_t$;
endow the latter set with some order and denote it
by ${\mathbf i} = (i_1,...,i_{m-1})$.
Let $L$ be as in Remark \ref{Rmk: degradation}; then
the restrictions $\mu_s|_L$ and $\mu_t|_L$ of $\mu_s$ and $\mu_t$
to $L$ are usual cross ratios of the variables
$p_j,p_k,p_s,p_t$ and $p_{j'},p_{k'},p_{s'},p_{t'}$,
respectively, and $(\mu_s|_L):(\mu_t|_L)$ is such
a cross ratio, as well. By Lemma \ref{Lm: simplices of cross ratios},
$\card(\{j,k,s,t\} \cap \{j',k',s',t'\}) =3$.
Consequently, $\card(\supp\mu_s \cap \supp\mu_t)
=\card {\supp_{\text{ess}}\mu_s} + 3 = m+2$.
Since $\card \supp\mu_s  = \card \supp\mu_t = m+3$,
it follows that $\supp\mu_s \not= \supp \mu_t$,
a contradiction.
Hence, for $s\ne t$ we have
$\card(\supp\limits_{\text{ess}}\mu_s
\cap\supp\limits_{\text{ess}}\mu_t) = m-2$.

Now let us show that $l\le m-1$.

Let ${\mathbf j}=(i_1,...,i_{m-2})=
\supp_{\text{ess}}\mu_1 \cap \supp_{\text{ess}}\mu_2$; then
$\mu_1 =e_{{\mathbf j},i;j,k,r,s}$ with certain $i,j,k,r,s$.
By Lemma \ref{cl: irreducibility} and
a straightforward computation, one can show that
$\mu_2\in D=\{e_{{\mathbf j},j;i,k,r,s},
e_{{\mathbf j},k;j,i,r,s},e_{{\mathbf j},r;j,k,i,s},
e_{{\mathbf j},s;j,k,r,i}\}$; the latter set
contains no pair of functions $\{\nu,\nu'\}$
that are vertices of the same simplex.

If $m=2$, then ${\mathbf j}=\varnothing$ and we must have
$l\le 1$, for otherwise it is easy to show that $\mu_3\in D$
and $\mu_2$ could not be a proper divisor of $\mu_3$.

Assume now that $m>2$.
Then $\card(\supp_{\text{ess}}\mu_1 \cap
\supp_{\text{ess}}\mu_2 \cap
\supp_{\text{ess}}\mu_t)<m-2$ for any $t>2$, since the equality
$\card(\supp_{\text{ess}}\mu_1 \cap
\supp_{\text{ess}}\mu_2 \cap
\supp_{\text{ess}}\mu_t)=m-2$ would imply
${\mathbf j}=\supp_{\text{ess}}\mu_1 \cap \supp_{\text{ess}}\mu_2=
\supp_{\text{ess}}\mu_1 \cap \supp_{\text{ess}}\mu_t$
and $\mu_t\in D$, which is impossible. Since the intersection of
any two of the three $m-1$ point sets
$\supp_{\text{ess}}\mu_1$, $\supp_{\text{ess}}\mu_2$ and
$\supp_{\text{ess}}\mu_t$ consists of $m-2$ points,
the intersection of all of them contains
at least $m-3$ points. Thus, $\card(\supp_{\text{ess}}\mu_1 \cap
\supp_{\text{ess}}\mu_2 \cap
\supp_{\text{ess}}\mu_t)=m-3$.

Furthermore, for any $t>2$, there is a unique
$i'\in\supp_{\text{ess}}\mu_1 \cap\supp_{\text{ess}}\mu_2$
such that $i'\not\in \supp_{\text{ess}}\mu_t$.
We shall show that $\mu_t$ is uniquely determined by this $i'$.

Let ${\mathbf i}=(i_1,...,i_{m-3})=
\supp_{\text{ess}}\mu_1 \cap
\supp_{\text{ess}}\mu_2 \cap \supp_{\text{ess}}\mu_t $
and $\mu_1 =e_{{\mathbf i},i,i';j,k,r,s}$.
By Lemma \ref{cl: irreducibility} and
a straightforward computation, one can show that
$\mu_2$ belongs to the set
$$
S=\{e_{{\mathbf i},i',j;i,k,r,s},
e_{{\mathbf i},i',k;j,i,r,s},e_{{\mathbf i},i',r;j,k,i,s},
e_{{\mathbf i},i',s;j,k,r,i}\}\,.
$$
Similarly,
$$
\mu_t\in T=\{e_{{\mathbf i},i,j;i',k,r,s},
e_{{\mathbf i},i,k;j,i',r,s},  e_{{\mathbf i},i,r;j,k,i',s},
e_{{\mathbf i},i,s;j,k,r,i'}\}\,.
$$
Neither $S$ nor $T$ contain
a pair of functions $\{\nu,\nu'\}$ such that $\nu|\nu'$.
Since
$$
\card(\supp\limits_{\text{ess}}\mu_2
\cap\supp\limits_{\text{ess}}\mu_t) = m-2\,,
$$
for every $\nu\in S$ there is only one $\nu'\in T$ such that
$\nu|\nu'$; this shows that $\mu_t$ is uniquely determined by
$i'$. It follows that $l=\dim\Delta\le m-1$.

Finally, in view of the above facts and the transitivity
of the ${\mathbf S}(n)$ action on $0$ dimensional simplices
(see Lemma \ref{Rmk: transitivity of cross ratios}),
the last statement of the part $(a)$ of the lemma is obvious.
\vskip7pt

\noindent $b)$
By Lemma \ref{Lm: m-generalized cross ratios lemma part a},
for $l<2$ the statement is obvious.
Suppose that $l \ge 2$. By  Lemma
\ref{Lm: m-generalized cross ratios lemma part a},
$\supp\limits_{\text{ess}}\mu_{s_0}
=\supp\limits_{\text{ess}}\mu_{t_0}$.
Assume that $\supp \mu_{\tilde s}=\supp \mu_{\tilde t}$
for some $\tilde s\ne\tilde t$.
Then, without loss of generality, we may assume
$\supp\mu_{\tilde s}=\supp\mu_{\tilde t}\ne\supp\mu_{s_0}$.
By  Lemma \ref{Lm: m-generalized cross ratios lemma part a},
this implies that
$\supp_{\text{ess}} \mu_{\tilde s}
=\supp_{\text{ess}} \mu_s $ and
$\supp_{\text{ess}} \mu_{\tilde t}
=\supp_{\text{ess}} \mu_s $; that is,
$\supp_{\text{ess}} \mu_{\tilde s}=
\supp_{\text{ess}} \mu_{\tilde t}$, which contradicts
part $(a)$ of the lemma. Thus,
$\supp_{\text{ess}}\mu_1 = ...= \supp_{\text{ess}}\mu_{l+1}$.
\vskip5pt

\noindent Let us prove that the simplex $\Delta$ is of the first
type, that is, there is a permutation
$\sigma$ such that $\{\sigma\mu_1,...,\sigma\mu_{l+1}\}
=\nabla_1^l$.
Set ${\mathbf i}
= \supp_{\text{ess}}\mu_1 = ...= \supp_{\text{ess}}\mu_{l+1}$,
$\mu_1=e_{{\mathbf i};j,k,r,s}$. By Lemma
\ref{Lm: m-generalized cross ratios lemma part a},
$\card(\supp \mu_1 \cap \supp \mu_2)=m+2$.
There is a unique index $t\in\supp\mu_2$ such that
$t\not\in\supp\mu_1$. Since $\mu_1|\mu_2$, Lemma
\ref{Lm: m-generalized cross ratios lemma part a} shows
that
$\mu_2\in D=\{e_{{\mathbf i};j,k,r,t},e_{{\mathbf i};k,j,s,t},
e_{{\mathbf i};r,s,j,t}, e_{{\mathbf i};s,r,k,t}\}$.
After a certain Kleinian permutation of the four last indices
$j,k,r,s$ in $\mu_1$, which never changes such a function,
and an appropriate renaming the indices in both $\mu_1$ and $\mu_2$,
we may assume that $\mu_1=e_{{\mathbf i};j,k,r,s}$
and $\mu_2=e_{{\mathbf i};j,k,r,t}$.

Let $p>2$; let us prove that $\supp\mu_1 \cap\supp \mu_2
=\supp\mu_1\cap\supp\mu_p=\supp\mu_2\cap\supp \mu_p$.
First, we shall prove that the index $t$ introduced above is not in
$\supp\mu_p$.
Suppose, on the contrary, that $t\in \supp \mu_p$.
Since $\mu_1|\mu_p$, by Lemma
\ref{Lm: m-generalized cross ratios lemma part a},
$\mu_p\in D$. But the set $D$ contains no pair of
determinant cross ratios which are proper divisors of each other,
which contradicts the condition $\mu_2|\mu_p$.
By the same reason, $s\not\in \supp \mu_p$.
Simple combinatorics show that
$\supp\mu_1 \cap\supp \mu_2
=\supp\mu_1\cap\supp\mu_p=\supp\mu_2\cap\supp \mu_p$.
Furthermore, there is a unique index $t_p\in\supp\mu_p$
such that $t_p\not\in\supp\mu_1$.
Since $\mu_1|\mu_p$, by Lemma
\ref{Lm: m-generalized cross ratios lemma part a},
$\mu_p\in \{e_{{\mathbf i};j,k,r,t_p},e_{{\mathbf i};k,j,s,t_p},
e_{{\mathbf i};r,s,j,t_p}, e_{{\mathbf i};s,r,k,t_p}\}$.
Since $\supp\mu_1\cap\supp\mu_2
=\supp\mu_1\cap\supp\mu_p$, it follows that
$\mu_p=e_{{\mathbf i};j,k,r,t_p}$.
The group $\mathbf{S}(n)$ is
$n$ times transitive on the set $\{1,...,n\}$,
thus, in view of Lemma
\ref{Rmk: transitivity of cross ratios}, there exists
a permutation $\sigma$ such that
$\{\sigma\mu_1,...,\sigma\mu_{l+1}\}=\nabla_1^l$.
This completes the proof of part $(b)$.
\vskip7pt

\noindent $c)$
The action of the permutation group preserves dimension
of simplices. By $(a)$ and $(b)$,
the ${\mathbf S}(n)$ orbit of any simplex contains a normal simplex.
We know that
$\dim \nabla_1^t \le n-m-3$ and
$\dim \nabla_2^t \le m-1$ (see
Definition \ref{def: normal simplices}).
That is, $\dim  L_{\vartriangle}({\mathcal E}^n(X))
= \max \{n-(m+3),m-1\}$.
\end{proof}
\vskip7pt

\noindent Notice that for $n>m+3$ the maximal possible
dimension of the simplices of the first type is $n-m-3\ge 1$;
in particular, the normal simplex of the first type
$\nabla_1^{n-m-3}$ is maximal in the sense of dimension.
Similarly, the maximal possible dimension of a simplex of
the second type is $m-1$ and the normal simplex of the second
type $\nabla_2^{m-1}$ is maximal.
In the following two lemmas we describe the stabilizers of
the maximal normal simplices $\nabla_1^{n-m-3}$ and
$\nabla_2^{m-1}$ in the group ${\mathbf S}(n)$ acting on
the set of all simplices of a given type and dimension;
since the latter action is transitive, these lemmas supply us
with an important information about the stabilizer of any maximal
simplex.

\begin{Lemma}\label{Lm: stabilizer of nabla1 maxdim}
Let $n>m+3$. The stabilizer ${\rm St}_{{\mathbf S}(n)}(\nabla_1^{n-m-3})$
of the ordered simplex $\nabla_1^{n-m-3}$ in the group ${\mathbf S}(n)$
conincides with the subgroup
${\mathbf S}(m-1)={\mathbf S}(\{1,...,m-1\})\subset{\mathbf S}(n)$.
\end{Lemma}

\begin{proof}
For any $s=m+3,...,n$, denote $I_s=(1,...,m+2,s)$ so that
$\nabla_1^{n-m-3} = \{e_{I_s}\}_{s=m+3}^n$.
Clearly ${\mathbf S}(m-1)
\subset{\rm St}_{{\mathbf S}(n)}(\nabla_1^{n-m-3})$.
Let $\sigma\in{\rm St}_{{\mathbf S}(n)}(\nabla_1^{n-m-3})$.
Then
$$
\aligned
e_{\{1,...,m-1\};m,m+1,m+2,{m+3}}&=
e_{I_{m+3}}=\sigma e_{I_{m+3}}
=\sigma  e_{\{1,...,m-1\};m,m+1,m+2,{m+3}}\\
&=e_{\{\sigma (1),...,\sigma (m-1)\};\sigma (m),\sigma (m+1),
\sigma (m+2),\sigma ({m+3})}
\endaligned
$$
and hence $\sigma$ is a disjoint product
$\sigma=\phi\psi\theta$,
where $\phi\in {\mathbf S}(m-1)$,
$\psi$ is
one of the four Kleinian permutations $\Id$, $(m,m+1)(m+2,m+3)$,
$(m,m+2)(m+1,m+3)$, $(m,m+3)(m+1,m+2)$
and $\theta\in{\mathbf S}(\{m+4,...,n\})$
(compare to Remark \ref{Rmk: when two DCR coincide}).
For any $t>m+3$
$$
\aligned
e_{\{1,...,m-1\};m,m+1,m+2,t}=
e_{I_t}&=\sigma e_{I_t}=\phi\psi\theta e_{I_t}
=\phi\psi\theta e_{\{1,...,m-1\};m,m+1,m+2,t}\\
&=e_{\{\phi (1),...,\phi (m-1)\};\psi (m),\psi (m+1),
\psi (m+2),\theta (t)}\\
&=e_{\{1,...,m-1\};\psi (m),\psi (m+1),
\psi (m+2),\theta (t)}\,.
\endaligned
$$
Consequently, $\psi =\Id$ and $\theta (t)=t$; since the latter
is true for any $t>m+3$, we see that $\theta=\Id$ and
$\sigma=\phi\in{\mathbf S}(m-1)$.
\end{proof}

\begin{Lemma}\label{Lm: stabilizer of nabla2 maxdim}
The stabilizer ${\rm St}_{{\mathbf S}(n)}(\nabla_2^{m-1})$
of the ordered simplex $\nabla_2^{m-1}$ in the group
${\mathbf S}(n)$ coincides with the subgroup
${\mathbf S}(\{m+4,...,n\})\subset{\mathbf S}(n)$.
\end{Lemma}

\begin{proof}
Of course, any element of ${\mathbf S}(\{m+4,...,n\})$
does not change $\nabla_2^{m-1}$. Let
$\sigma\in{\rm St}_{{\mathbf S}(n)}(\nabla_2^{m-1})$.
Denote $I=(1,...,m+3)$; then
$\nabla_2^{m-1}=\{(i,m)e_{I}\}_{i=1}^m$, where, as usual,
$(i,t)$ denotes the transposition of two indices $i,t$;
furthermore,
$(i,m)e_{I}=\sigma(i,m)e_{I}$, i.e., $e_{I}=(i,m)\sigma(i,m)e_{I}$.
For $i=m$ this means that
$$
\aligned
e_{\{1,...,m-1\};m,m+1,m+2,m+3}&=e_{I}=\sigma e_{I}
=\sigma e_{\{1,...,m-1\};m,m+1,m+2,m+3}\\
&=e_{\{\sigma (1),...,\sigma (m-1)\};\sigma (m),\sigma (m+1),
\sigma (m+2),\sigma ({m+3})}
\endaligned
$$
and hence $\sigma$ is a disjoint product
$\sigma=\theta\phi\psi$, where $\theta\in {\mathbf S}(m-1)$,
$\phi$ is one of the four Kleinian permutations
$\Id$, $(m,m+1)(m+2,{m+3})$,
$(m,m+2)(m+1,m+3)$, $(m,{m+3})(m+1,m+2)$ and
$\psi\in{\mathbf S}(\{m+4,...,n\})$
(compare to Remark \ref{Rmk: when two DCR coincide}).
For any $i=1,...,m-1$, we have
$$
\aligned
e_{\{1,...,i-1,m,i+1,...,m-1\};i,m+1,m+2,m+3}
&=(i,m)e_{\{1,...,m-1\};m,m+1,m+2,m+3}\\
&\hskip-89pt=(i,m)e_{I}
                    =\sigma(i,m)e_{I}=\theta\phi\psi(i,m) e_{I}\\
&\hskip-89pt=\theta\phi\psi(i,m)e_{\{1,...,m-1\};m,m+1,m+2,m+3}\\
&\hskip-89pt=e_{\{\theta (1),...,\theta (i-1),\phi (m),
\theta (i+1),...,\theta (m-1)\};\theta (i),
\phi (m+1),\phi (m+2),\phi (m+3)}\,.
\endaligned
$$
The latter equality is fulfilled if and only if
$\phi=\Id$ and $\theta(i)=i$; since the latter
is true for any $i=1,...,m-1$, we see that $\theta=\Id$ and
$\sigma=\psi\in{\mathbf S}(\{m+4,...,n\})$.
\end{proof}

\subsection{Maps of $L_{\vartriangle}({\mathcal E}^n(X,gp))$
induced by holomorphic endomorphisms}
\label{ss: Induced maps on the complex L(En)}
Here we prove that any strictly equivariant holomorphic
map $f\colon{\mathcal E}^n(X,gp)\to{\mathcal E}^n(X,gp)$
induces a simplicial map
$f^*\colon L_{\vartriangle}({\mathcal E}^n(X,gp))\to
L_{\vartriangle}({\mathcal E}^n(X,gp))$.
Then we show that, for sufficiently large $n$, the
vertices of certain $m$ simplices of maximal dimension
are, up to a permutation of $q_1,..,q_n$, fixed points
of the map $f^*$.
\vskip7pt

\noindent We need the following result similar
to Lemma 6.1 and Corollary 6.2 from \cite{LinSphere}.

\begin{Lemma}\label{Lm: the image of ratio is ratio}
Let $X$ be either $\Bbb{CP}^m$ or $\Bbb{C}^m$.
Let $f\colon\Cal{E}^n(X,gp)\rightarrow \Cal{E}^n(X,gp)$ be
a strictly equiva\-riant ho\-lo\-morphic map.
If $\lambda\in \hskip-1pt L({\mathcal E}^n(X,gp))$ then
$ \lambda\circ f \in \hskip-1pt L({\mathcal E}^n(X,gp))$.
\end{Lemma}

\begin{proof}
It suffices to show that $\lambda\circ f\ne\const$
for any $\lambda\in L({\mathcal E}^n(X,gp))$.
Assume to the contrary that
$\lambda\circ f =c\in \Bbb{C}\setminus\{0,1\}$
for some $\lambda\in L_{\vartriangle}({\mathcal E}^n(X,gp))$.
Then $(\lambda\circ f)(\sigma q)\equiv c$
for all $\sigma \in{\mathbf S}(n)$;
as $f$ is strictly equivariant,
this implies $\lambda(\sigma f(q))\equiv c$ for each
$\sigma\in{\mathbf S}(n)$.
Since $\lambda\in L({\mathcal E}^n(X,gp))$, it follows that
$\lambda^{-1},1-\lambda \in L({\mathcal E}^n(X,gp))$.
By Lemma \ref{Rmk: transitivity of cross ratios},
there are $s_1,s_2\in{\mathbf S}(n)$
such that $\lambda^{-1}=s_1\lambda$ and
$1-\lambda=s_2\lambda$.
Set $\sigma_i=\alpha^{-1}(s_i^{-1})$,
where $\alpha$ is an automorphism of
${\mathbf S}(n)$ related to the strictly equivariant
mapping $f$ and $i=1,2$.
Clearly,
$c^{-1}\equiv\lambda^{-1}(f(q))\equiv
s_1\lambda(f(q))\equiv\lambda(s_1^{-1}f(q))
\equiv\lambda(f(\sigma_1 q))\equiv c$ and
$1-c\equiv(1-\lambda)(f(q))\equiv
s_2\lambda(f(q))\equiv\lambda(s_2^{-1}f(q))\equiv
\lambda(f(\sigma_2 q))\equiv c$.
Both of the above equalities can not occur
at the same time, this is a contradiction.
\end{proof}

\begin{Corollary}\label{Crl: induced map}
Let $X=\Bbb{CP}^m$ or $\Bbb{C}^m$ and $n>m+2$. Let
$ f\colon\Cal{E}^n(X,gp)\rightarrow \Cal{E}^n(X,gp)$
be a strictly equivariant holomorphic map.
Then $f$ induces a simplicial mapping $f^*$
whose restriction to each simplex
$\Delta\subseteq L(\Cal{E}^n(X,gp))$ is injective, and
dimension of a simplex does not change under this
transformation.
\end{Corollary}

\begin{proof}
By Lemma \ref{Lm: the image of ratio is ratio},
For any $\lambda\in L({\mathcal E}^n(X,gp))$ we have
$\lambda\circ f \in L({\mathcal E}^n(X,gp))$. Thus
$f$ induces a map
$f^*\colon L(\Cal{E}^n(X,gp))\to L(\Cal{E}^n(X,gp))$ defined by
$$
L({\mathcal E}^n(X,gp))\ni\lambda
\mapsto f^*(\lambda ) = \lambda\circ f\in
L({\mathcal E}^n(X,gp))\,.
$$
If $\lambda,\mu,\nu\in L({\mathcal E}^n(X,gp))$ and $\lambda=\mu/\nu$
then
$$
f^*(\mu)/f^*(\nu) = f^*(\mu/\nu) =
f^*(\lambda)\in L({\mathcal E}^n(X,gp))\,.
$$
Therefore,
$f^*$ is a simplicial map, its restriction to each simplex
$\Delta \in  L_{\vartriangle}(\Cal{E}^n(X,gp))$ is injective and
dimension of a simplex does not change under this
transformation.
\end{proof}

\begin{Lemma}\label{Lm: f* isomorphism}
Let  $X$ be either $\Bbb{CP}^m$ or $\Bbb{C}^m$,
$n\ge m+3$ and $n\ne2m+2$.
Let $ f\colon\Cal{E}^n(X,gp)\rightarrow \Cal{E}^n(X,gp)$ be
a strictly equivariant holomorphic map.
Then the induced simplicial map $f^*$ is an automorphism of
the complex $L_\vartriangle(\Cal{E}^n)$. Moreover,
$f^*$ preserves the type of simplices.
\end{Lemma}

\begin{proof}
The set $L(\Cal{E}^n)$ of all vertices of
$L_\vartriangle(\Cal{E}^n)$ is finite. Hence, to prove that
the simplicial map $f^*$ is an automorphism of
the complex $L_\vartriangle(\Cal{E}^n)$, it suffices to show that
the map $f^*\colon L(\Cal{E}^n)\to L(\Cal{E}^n)$ is
bijective, which, by the same finiteness reason,
is equivalent to its surjectivity.

Let $\alpha\in\Aut({\mathbf S}(n))$ be the automorphism
related to our strictly equivariant holomorphic map $f$ so that
$f(\sigma q)= \alpha(\sigma)f(q)$ and $f^*(\sigma\mu)=
\alpha^{-1}(\sigma)[f^*(\mu)]$ for all $\sigma\in{\mathbf S}(n)$
and $\mu\in L(\Cal{E}^n)$.

Let $\mu\in L(\Cal{E}^n)$ and $\nu=f^*(\mu)$. By Lemma
\ref{Rmk: transitivity of cross ratios}, there exists
$\sigma\in {\mathbf S}(n)$ such that $\sigma \nu =\mu$.
Set $\lambda= \alpha(\sigma)\mu\in  L(\Cal{E}^n)$;
then $\mu=\sigma\nu=\sigma (f^*(\mu))=f^*(\alpha(\sigma)\mu)
=f^*(\lambda)$, which proves that the map
$f^*\colon L(\Cal{E}^n)\to L(\Cal{E}^n)$ is surjective and
thereby bijective. It follows that $f^*$ is a simplicial
automorphism of the finite complex $L_\vartriangle(\Cal{E}^n)$.
\vskip5pt

\noindent Let us prove now that $f^*$ preserves the type of
simplices. First assume that $n>2m+2$.
We start with the normal simplex
$\nabla_1^{n-m-3}$ and its faces. Since $f^*$ preserves dimension,
$\dim f^*(\nabla_1^{n-m-3})=n-m-3>m-1$ and
Lemma \ref{Lm: m-generalized cross ratios lemma}$(a,b)$
shows that the simplex $f^*(\nabla_1^{n-m-3})$ is of the first
type. Any normal simplex of the first type
$\nabla_1^{l}$ is a face of $\nabla_1^{n-m-3}$ and
any face of a simplex of the first type is also
a simplex of the first type (see Definition
\ref{def: normal simplices}).
Since
$f^*(\nabla_1^{l})$ is a face of the simplex
$f^*(\nabla_1^{n-m-3})$ which is of the first type,
$f^*(\nabla_1^{l})$ is of the first type.

Now let $\Delta\in L_\vartriangle(\Cal{E}^n)$ be
any $l$-simplex of the first type. It follows from Definition
\ref{def: normal simplices}
that there is $\sigma\in {\mathbf S}(n)$ such that
$\sigma \nabla_1^{l}=\Delta$. Therefore
the simplex $f^*(\Delta)=f^*(\sigma \nabla_1^{l})=
\alpha^{-1}(\sigma) f^*(\nabla_1^{l})$
is of the first type. Thus, $f^*$ carries simplices of the first
type to simplices of the first type.

Since the simplicial map $f^*$ is an automorphism
of the finite complex $L_\vartriangle(\Cal{E}^n)$,
$f^*$ is bijective on the set of all simplices
of positive dimension.
By Lemma \ref{Lm: m-generalized cross ratios lemma}$(a,b)$ and
Definition \ref{def: normal simplices}, the latter set is a
disjoint union of two its subsets consisting of all simplices
of the first and the second type, respectively. Therefore, it
follows from what was proved above that $f^*$ carries
simplices of the second type to simplices of the second type.
This completes the proof of the case $n>2m+2$
\vskip5pt

\noindent When $n<2m+2$, we consider the normal simplex
$\nabla_2^{m-1}$ and its faces. Now we have
$\dim f^*(\nabla_2^{m-1})=m-1>n-m-3$ and, by 
Lemma \ref{Lm: m-generalized cross ratios lemma}$(a,b)$,
the simplex $f^*(\nabla_2^{m-1})$ is of the second
type. Any normal simplex of the second type
$\nabla_2^{l}$ is a face of $\nabla_2^{m-1}$ and
any face of a simplex of the second type is also
a simplex of the second type (see Definition
\ref{def: normal simplices}).
Since $f^*(\nabla_2^{l})$ is a face of the simplex
$f^*(\nabla_2^{m-1})$ of the second type,
$f^*(\nabla_2^{l})$ is also of the second type. Thus,
$f^*$ preserves the type of the simplices of the second type.

Any $l$-simplex $\Delta\in L_\vartriangle(\Cal{E}^n)$ 
of the second type may be represented as $\Delta=\sigma \nabla_2^{l}$
with some $\sigma\in {\mathbf S}(n)$. Therefore
the simplex $f^*(\Delta)=f^*(\sigma \nabla_2^{l})=
\alpha^{-1}(\sigma) f^*(\nabla_2^{l})$
is of the second type. Precisely as above,
using the bijectivity of $f^*$ on the set of all simplices
of positive dimension and taking into account
Lemma \ref{Lm: m-generalized cross ratios lemma}$(a,b)$
and Definition \ref{def: normal simplices}, we come to conclusion
that $f^*$ carries simplices of the first type to simplices of
the first type, which completes the proof of the lemma. 
\end{proof}

\begin{Remark}\label{Rmk: f* and S(n) action}
Let $\alpha\in\Aut{\mathbf S}(n)$ be the
automorphism related to a strictly equivariant
holomorphic endomorphism $f$ of ${\mathcal E}^n$.
It is easily seen that for any $\rho\in{\mathbf S}(n)$ and any
function $\lambda$ on ${\mathcal E}^n$ we have
\begin{equation}\label{eq: action of permutation on composition}
(\alpha(\rho)f)^*(\lambda)=\lambda\circ [\alpha(\rho)f)]
=\lambda\circ f\circ\rho=\rho^{-1}[\lambda\circ f]
=\rho^{-1}[f^*(\lambda)]\,,
\end{equation}
where $\rho^{-1}[\lambda\circ f]$ is the result of the action
of the permutation $\rho^{-1}$ on the function
$\lambda\circ f=f^*(\lambda)$.
Changing $\rho$ with $\rho^{-1}$, we obtain
\begin{equation}
\label{eq: action of permutation on composition additional}
\rho[f^*(\lambda)]
=(\alpha(\rho^{-1})f)^*(\lambda)=
\lambda\circ [\alpha(\rho^{-1})f)]
=(\alpha(\rho)\lambda)\circ f =f^*(\alpha(\rho)\lambda)\,.
\end{equation}
\hfill $\bigcirc$
\end{Remark}

\begin{Remark}
Let $f=f_{\tau,\sigma}\colon\Cal{E}^n\rightarrow \Cal{E}^n$ be
a strictly equivariant tame holomorphic map and let
$\mu\in L(\Cal{E}^n)$. Since the function $\mu$ is
$\mathbf{PSL}(m+1,\mathbb{C})$ invariant,
we have
$(f^*\mu)(q)=\mu(f(q))=\mu(\sigma\tau(q)q)=
\mu(\tau(q)\sigma q)=\mu(\sigma q)$. Taking into account
that $f$ is strictly equivariant, we obtain
\begin{equation}\label{eq: f^* acts on L as permutation}
\aligned
\mu(q)=(f^*\mu)(\sigma^{-1}q)
=\mu(f(\sigma^{-1}q))=\mu(\alpha(\sigma^{-1})f(q))
=[(\alpha(\sigma^{-1})f)^*\mu](q)\,.
\endaligned
\end{equation}
The latter formula shows that all the vertices
$\mu\in L(\Cal{E}^n)$ of the complex $L_\vartriangle(\Cal{E}^n)$
are fixed points of the map $(\alpha(\sigma^{-1})f)^*$;
in other words, {\em the action of the mapping $f^*$ itself
on the set $L(\Cal{E}^n)$ of all vertices of the complex
$L_\vartriangle(\Cal{E}^n)$ coincides with the action
of a certain permutation}.

In the theorem below we prove a weaker property
of all strictly equivariant holomorphic endomorphisms $f$
of $\Cal{E}^n$. It says that for any such $f$
there is a permutation $\rho\in{\mathbf S}(n)$
such that the vertices of certain special
simplices of $L_\vartriangle(\Cal{E}^n)$ are fixed points
of the map $(\rho f)^*$; that is, the action of $f^*$
on those particular vertices of $L_\vartriangle(\Cal{E}^n)$
coincides with the action of a certain permutation.
This prelimenary result will play important part in the proof
of Theorem \ref{Thm: endomorphisms of En(X)}.
\hfill $\bigcirc$
\end{Remark}

\noindent According to Notation \ref{Nt: first m members},
${\mathbf m}(\hat s) = (1,...,\hat{s}, ...,m)$ for any $s=1,..,m$;
sometimes we write $\widehat{\mathbf m}$
instead of ${\mathbf m}(\hat m)$. We use the standard notation
$(i,j)$ for the transposition of two distinct elements
$i,j\in\{1,...,n\}$. Notice that for $i=1,...,m$
the permutation $(i,m)\nabla_1^{n-m-3}$ of the simplex
$\nabla_1^{n-m-3}=\{e_{\widehat{\mathbf m};m,m+1,m+2,s}\}_{s=m+3}^n$
is the simplex $\{e_{{\mathbf m}(\hat i);i,m+1,m+2,s}\}_{s=m+3}^n$.

\begin{Theorem}
\label{Thm: cross ratio is invariant under automorphism for dim>2}
Let  $X$ be either $\Bbb{CP}^m$ or $\Bbb{C}^m$, $n\ge m+3$ and
$n\ne 2m+2$.
Let $ f\colon\Cal{E}^n(X,gp)\rightarrow \Cal{E}^n(X,gp)$ be
a strictly equivariant holomorphic map.
There exists a permutation
$\rho\in {\mathbf S}(n)$ such that
$(\rho f)^*(e_{{\mathbf m}(\hat r);r,m+1,m+2,s})
=e_{{\mathbf m}(\hat r);r,m+1,m+2,s}$
for any $r\in \{1,...,m\}$ and $s\in\{m+3,...,n\}$.
In other words,
the map $(\rho f)^*$ is identical on each of the simplices
\begin{equation}\label{eq: fixed simplices}
\{e_{{\mathbf m}(\hat 1);1,m+1,m+2,s}\}_{s=m+3}^n,
...,\{e_{{\mathbf m}(\hat m);m,m+1,m+2,s}\}_{s=m+3}^n\,.
\end{equation}
\end{Theorem}

\begin{proof}
First assume that $n>m+3$.
Notice that the last simplex in the above list, namely
$$
\{e_{{\mathbf m}(\hat m);m,m+1,m+2,s}\}_{s=m+3}^n
=\{e_{\widehat{\mathbf m};m,m+1,m+2,s}\}_{s=m+3}^n
=\nabla_1^{n-m-3}\,,
$$
is the normal $(n-m-3)$-simplex
of the first type (see Definition \ref{def: normal simplices}).
To simplify the notation, for any $s=m+3,...,n$,
let us set $I_s=(1,...,m+2,s)$;
notice that $I_s$ is the support of
the determinant cross ratio $e_{\widehat{\mathbf m};m,m+1,m+2,s}$,
which is one of the vertices of $\nabla_1^{n-m-3}$.

By Lemma \ref{Lm: f* isomorphism},
$f^*$ preserves the type of simplices; hence
$f^*(\nabla_1^{n-m-3})$ is a simplex of the first type and
there is a permutation $\theta$ such that
$(\theta f)^*(\nabla_1^{n-m-3})=\nabla_1^{n-m-3}$.
Clearly, it would be sufficient to prove the theorem
for the strictly equivariant map $\theta f$.
Therefore, without loss of generality, we may
from the very beginning assume that our original map $f$ satisfies
$f^*(\nabla_1^{n-m-3})=\nabla_1^{n-m-3}$.
Since we deal with ordered simplices, the latter relation
means that all vertices of the simplex $\nabla_1^{n-m-3}$ are
fixed points of the map $f^*$, that is,
$$
f^*(e_{\widehat{\mathbf m};m,m+1,m+2,s})
=e_{\widehat{\mathbf m};m,m+1,m+2,s}\,,
$$
or, which is the same,
\begin{equation}\label{eq: vertices of nabla1(n-m-3) are fixed points of f*}
f^*(e_{I_s})
=e_{I_s}\quad\text{\rm for all}\ \
s=m+3,...,n\,,
\end{equation}
Since $f$ is strictly equivariant, there is
$\alpha\in\Aut({\mathbf S}(n))$
such that $f(\sigma q)= \alpha(\sigma) f(q)$ for
every $\sigma \in {\bf S}(n)$.
Consequently, for $1\le i < m$ and $m \le t \le n$
we have
\begin{equation}\label{eq: the image of permuted simplex}
f^*((i,t)\nabla_1^{n-m-3})=
\alpha^{-1}((i,t))\nabla_1^{n-m-3}\,,
\end{equation}
where $(i,t)$ is the transposition of $i$ and $t$.

For any $s=m+3,...,n$, set $I_s=(1,...,m+2,s)$;
notice that $I_s$ is the support of
the determinant cross ratio $e_{\widehat{\mathbf m};m,m+1,m+2,s}$,
which is one of the vertices of $\nabla_1^{n-m-3}$.
\vskip5pt

\noindent The permuted simplices
$$
\aligned
\Delta_{s,1}&=(m+3,s)\nabla_2^{m-1} =
\{(1,m)e_{I_s},(2,m)e_{I_s},...,(m-1,m)e_{I_s},e_{I_s}\}\,,\\
\Delta_{s,2}&=(m+3,s)(m,m+1)(m+2,m+3)\nabla_2^{m-1}\\
&=\{(1,m+1)e_{I_s},(2,m+1)e_{I_s},...,(m-1,m+1)e_{I_s},e_{I_s}\}\,,\\
\Delta_{s,3}
&=(m+3,s)(m,m+2)(m+1,m+3)\nabla_2^{m-1}\\
&=\{(1,m+2)e_{I_s},(2,m+2)e_{I_s},...,(m-1,m+2)e_{I_s},e_{I_s}\}\,,\\
\Delta_{s,4}
&=(m+3,s)(m,m+3)(m+1,m+2)\nabla_2^{m-1}\\
&=\{(1,s)e_{I_s},(2,s)e_{I_s},...,(m-1,s)e_{I_s},e_{I_s}\}
\endaligned
$$
are the simplices of the second type.
By Lemma \ref{Lm: f* isomorphism},
all $f^*(\Delta_{s,\kappa})$,
$\kappa=1,2,3,4$, are simplices of the second type.
Since $f^*$ preserves vertices of
$\nabla_1^{n-m-3}$, we have $f^*(e_{I_s})=e_{I_s}$
and hence each of the simplices $f^*(\Delta_{s,\kappa})$,
$\kappa=1,2,3,4$, contains the vertex $e_{I_s}$ whose
essential support is $\widehat{\mathbf m}=(1,...,m-1)$.
By Lemma \ref{Lm: m-generalized cross ratios lemma}(a),
the essential supports
of all vertices of $f^*(\Delta_{s,\kappa})$ but $e_{I_s}$
are not equal to $\supp_\text{ess} e_{I_s}=\widehat{\mathbf m}$.


\vskip5pt

\noindent{\bf Claim.}
{\sl There exists $\sigma \in {\mathbf S}(m-1)
={\mathbf S}(\{1,...,m-1\})\subset {\mathbf S}(n)$
such that for any $s>m+2$ the couple of the simplices
$f^*(\Delta_{s,2})$, $f^*(\Delta_{s,3})$ coincides with
the couple of the simplices
$\sigma\Delta_{s,2}$, $\sigma \Delta_{s,3}$.}
\vskip5pt

\noindent{\em Proof of Claim.}
We divide the proof to four steps.
\vskip5pt

{\bf Step 1.} Pick some $s\ge m+3$. Let
$$
\aligned
\Delta_s&=f^*(\Delta_{s,2})\\
&=
\{f^*((1,m+1)e_{I_s}),f^*((2,m+1)e_{I_s}),...,
f^*((m-1,m+1)e_{I_s}),f^*(e_{I_s})\}\\
&= \{f^*((1,m+1)e_{I_s}),f^*((2,m+1)e_{I_s}),...,
f^*((m-1,m+1)e_{I_s}),e_{I_s}\}\,;
\endaligned
$$
then $\Delta_s$ is a simplex of the second type and
Lemma \ref{Lm: m-generalized cross ratios lemma}(a,b) shows
that
\begin{equation}\label{eq: equality of supports}
\supp f^*((1,{m+1})e_{I_s})=...=\supp f^*((m-1,{m+1})e_{I_s})
=\supp e_{I_s}=I_s\,.
\end{equation}
The simplex $\Delta_{s,2}$ is also of the second type
and, according to Definition \ref{def: normal simplices},
it may be carried to the simplex $\Delta_s$
by a permutation $\phi_s\in{\mathbf S}(n)$.
The vertices of both these simplices depend
only on the vector variables $q_1,...,q_{m+2},q_s$
(see Section \ref{ss: Notation and definitions} and
Definition \ref{Def: determinant cross ratio});
by Remark
\ref{rm: simplices involving the same set of vector variables},
$\phi_s$ may be chosen in the subgroup ${\mathbf S}(I_s)
={\mathbf S}(\{1,...,m+2,s\})\subset{\mathbf S}(n)$ and,
by Lemma \ref{Lm: stabilizer of nabla2 maxdim},
such a permutation is unique.
In particular, we have
$$
\aligned
e_{I_s}&=\phi_s e_{I_s}
=\phi_s e_{\{1,...,m-1\};m,m+1,m+2,s}\\
&=e_{\{\phi_s(1),...,\phi_s(m-1)\};\phi_s(m),\phi_s(m+1),
\phi_s(m+2),\phi_s(s)}
\endaligned
$$
and hence $\phi_s$ is a disjoint product $\phi_s=\sigma_s\theta_s$,
where $\sigma_s\in {\mathbf S}(m-1)$ and $\theta_s$ is
one of the four Kleinian permutations $\Id$, $(m,m+1)(m+2,s)$,
$(m,m+2)(m+1,s)$, $(m,s)(m+1,m+2)$
(compare to Remark \ref{Rmk: when two DCR coincide}).
Consequently, for any $i=1,...,m-1$ we have
\begin{equation}
\label{eq: the action of Kleinian permutation on permuted cross ratio}
\aligned
\theta_s (i,m+1)e_{I_s}
&=\theta_s (i,m+1)e_{\{1,...,m-1\};m,m+1,m+2,s}\\
&=\theta_s e_{\{1,...,i-1,m+1,i+1,...,m-1\};m,i,m+2,s} \\
&=e_{\{1,...,i-1,\theta_s(m+1),i+1,...,m-1\};
\theta_s(m),i,\theta_s(m+2),\theta_s(s)}\,;
\endaligned
\end{equation}
the latter function, in turn, must be
one of the determinant cross ratios
$$
(i,{m+1})e_{I_s}\,,  \  \ (i,{m})e_{I_s}\,, \ \
(i,s)e_{I_s}\,, \  \ (i,{m+2})e_{I_s}\,.
$$
This means that $\theta_s\Delta_{s,2} = \Delta_{s,j_s}$
for a certain $j_s\in\{1,2,3,4\}$ and hence
\begin{equation}\label{eq: intermediate formula for f*(Delta(s,2)}
\Delta_s=f^*(\Delta_{s,2})=\sigma_s\theta_s\Delta_{s,2}
=\sigma_s \Delta_{s,j_s}\,.
\end{equation}
Since $\sigma_s\in{\mathbf S}(m-1)$ does not touch the indices
$m,m+1,...,n$, relation
\eqref{eq: intermediate formula for f*(Delta(s,2)}
determines both $j_s$ and $\theta_s$;
moreover, $j_s$ is uniquely determined by the value
of $\theta_s$ on any one of the numbers $m,m+1,m+2,s$.
In particular, we see that
the value $\theta_s(m+1)\in\{m,m+1,m+2,s\}$
determines $j_s$; more precisely,
\vskip5pt

$(i)$ the ordered couple $(j_s;\theta_s(m+1))$ is one
of the four ordered couples
\begin{equation}\label{eq: four ordered couples}
(1;m)\,,\ (2;m+1)\,, \ (3;m+2)\,, \ (4,s)\,.
\end{equation}
\vskip5pt

{\bf Step 2.} Let us show now that both the permutation
$\sigma_s\in{\mathbf S}(m-1)$ and the index $j_s\in\{1,2,3,4\}$
in \eqref{eq: intermediate formula for f*(Delta(s,2)}
do not depend on $s$ and, moreover, $j_s\ne 4$.

We start with $\sigma_s$. Fix some $i\in\{1,...,m-1\}$.
According to \eqref{eq: the image of permuted simplex},
the simplex $f^*((i,m+1)\nabla_1^{n-m-3})
=\alpha^{-1}((i,m+1))\nabla_1^{n-m-3}$ is of the first type;
the functions
\begin{equation*}
f^*((i,{m+1}) e_{I_{m+3}}),f^*((i,{m+1}) e_{I_{m+4}}),...,
f^*((i,{m+1}) e_{I_{n}})
\end{equation*}
are its vertices, and Lemma \ref{Lm: m-generalized cross ratios lemma}(a,b)
shows that
\begin{equation}
\label{eq: equality of essential support of permuted simplex of the first type}
\supp_{ess} f^*((i,{m+1}) e_{I_{m+3}}) = ...
=\supp_{ess} f^*((i,{m+1}) e_{I_n})\,.
\end{equation}
Furthermore, the function $f^*((i,{m+1}) e_{I_{s}})$
is the $i^{\text{\rm th}}$ vertex of the simplex $f^*(\Delta_{s,2})$;
by \eqref{eq: intermediate formula for f*(Delta(s,2)},
$f^*(\Delta_{s,2})=\sigma_s\theta_s\Delta_{s,2}$ and hence, by
\eqref{eq: the action of Kleinian permutation on permuted cross ratio},
$$
\aligned
f^*((i,{m+1})& e_{I_{s}})= \sigma_s\theta_s(i,{m+1}) e_{I_{s}}\\
&=\sigma_s e_{\{1,...,i-1,\theta_s(m+1),i+1,...,m-1\};
\theta_s(m),i,\theta_s(m+2),\theta_s(s)}\\
&=e_{\{ \sigma_s(1),..., \sigma_s(i-1),\theta_s(m+1), \sigma_s(i+1),
..., \sigma_s(m-1)\};\theta_s(m),\sigma_s(i),\theta_s(m+2),
\theta_s(s)}\,.
\endaligned
$$
Therefore, by
\eqref{eq: equality of essential support of permuted simplex of the first type},
the set
\begin{equation}\label{eq: essential support 1}
\aligned
\Sigma&\Def \supp_{ess}f^*((i,{m+1})e_{I_{s}})\\
&=\{\sigma_{s}(1),...,\sigma_{s}(i-1),\theta_s(m+1),
\sigma_{s}(i+1),...,\sigma_{s}(m-1)\}
\endaligned
\end{equation}
does not depend on $s$. The only element of
$\Sigma$ that is not in $\{1,...,m-1\}$ is
$\theta_{s}(m+1)$; hence $\theta_{s}(m+1)$ does not depend on $s$
and the same is true for the set $\Sigma'\Def
\{\sigma_{s}(1),...,\widehat{\sigma_{s}(i)},...,\sigma_{s}(m-1)\}$
of all elements of $\Sigma$ but $\theta_{s}(m+1)$.
In fact the set $\Sigma'$ consists of all numbers $1,...,m-1$
but $\sigma_s(i)$; thus, $\sigma_s(i)$ also does not depend on $s$.
Since this is the case for any $i\in\{1,...,m-1\}$, the
the permutation $\sigma\Def\sigma_s\in{\mathbf S}(m-1)$
does not depend on $s$.

Now we turn to the index $j_s\in\{1,2,3,4\}$.
$\theta_s$ is a Kleinian permutation of
$m,m+1,m+2,s$ and we have already proved that the element
$$
\theta_s(m+1)\in\{m,m+1,m+2,s\}
$$
does not depend on $s\in\{m+3,...,n\}$; thus,
$$
\theta_{m+3}(m+1)=...=\theta_n(m+1)
\in \bigcap_{s=m+3}^n \{m,m+1,m+2,s\}=\{m,m+1,m+2\}
$$
and hence $\theta_s(m+1)\ne s$.
According to $(i)$, this means that $j_s\ne 4$
and $j_s$ does not depend on $s$.

Thus, as a result of steps 1 and 2, we know that
\vskip5pt

$(ii)$ there are a permutation $\sigma\in {\mathbf S}(m-1)$
and $j\in\{1,2,3\}$ such that
\begin{equation}\label{eq: how it looks like1}
f^*(\Delta_{s,2}) = \sigma \Delta_{s,j}\quad
\text{\rm for any} \ \ s\in\{m+3,...,n\}\,.
\end{equation}
\vskip5pt

{\bf Step 3.} In a similar way, one can show that
\vskip5pt

$(ii')$ there are a permutation $\sigma'\in {\mathbf S}(m-1)$
and $j'\in\{1,2,3\}$ such that
\begin{equation}\label{eq: how it looks like1'}
f^*(\Delta_{s,3}) = \sigma' \Delta_{s,j'}\quad
\text{\rm for any} \ \ s\in\{m+3,...,n\}\,.
\end{equation}
\vskip5pt

{\bf Step 4.} Let us prove now that $\sigma=\sigma'$,
$j,j'\ne 1$ and $j\ne j'$;
this will complete the proof of Claim.

By Lemma \ref{Lm: 3 facts about strong equivariance multi dim}(d),
for any $s\in \{m+3,...n\}$ the function
$e_{I_s}$ admits the following two representations as products
of two determinant cross ratios:
\begin{equation}\label{eq: the prefunction relation}
\aligned
e_{I_s}=
((i,{m}) e_{I_s})\cdot
((i,s) e_{I_s})=
((i,{m+1}) e_{I_s})\cdot
((i,{m+2}) e_{I_s})\,.
\endaligned
\end{equation}
Hence, by \eqref{eq: vertices of nabla1(n-m-3) are fixed points of f*}
\begin{equation}\label{eq: the function relation}
\aligned
e_{I_s}&=f^*(e_{I_s})=f^*((i,{m}) e_{I_s})f^*((i,s) e_{I_s})\\
&=f^*((i,{m+1}) e_{I_s})f^*((i,{m+2}) e_{I_s})\,.
\endaligned
\end{equation}

Suppose that $j=1$. Then, according to
\eqref{eq: how it looks like1}, $f^*((i,{m+1})e_{I_s})$ is the
$i^{\rm th}$ vertex of the simplex
$f^*(\Delta_{s,2})=\sigma \Delta_{s,1}$
and hence $f^*((i,{m+1})e_{I_s})=(\sigma(i),m)e_{I_s}$.
By \eqref{eq: the function relation} and
\eqref{eq: the prefunction relation}, we have
\begin{equation}\label{eq: the its vertex of f*(Delta(s,3))}
f^*((i,m+2) e_{I_s})=
\frac{e_{I_s}}{f^*((i,{m+1}) e_{I_s})} =
\frac{e_{I_s}}{(\sigma(i),m) e_{I_s}}
= (\sigma(i),s) e_{I_s}\,.
\end{equation}
By \eqref{eq: how it looks like1'} and
\eqref{eq: the its vertex of f*(Delta(s,3))},
the function $(\sigma(i),s) e_{I_s}$ must be the
$i^{\rm th}$ vertex of the simplex
$f^*(\Delta_{s,3})=\sigma' \Delta_{s,j'}$; the latter
shows that $j'=4$, which contradicts $(ii')$. Thus $j\ne 1$.

The proof of the inequality $j'\ne 1$
is similar to the above one. Consequently,
$j,j'\in\{2,3\}$.

Now we turn to the permutations $\sigma$ and $\sigma'$.

Suppose first that $j=2$. Then, according to
\eqref{eq: how it looks like1}, each $f^*((i,m+1)e_{I_s})$,
$i=1,...,m-1$, is a vertex of the simplex
$f^*(\Delta_{s,2})=\sigma \Delta_{s,2}$;
hence $f^*((i,m+1)e_{I_s})=(\sigma(i),m+1)e_{I_s}$.
By \eqref{eq: the function relation} and
\eqref{eq: the prefunction relation}, we have
\begin{equation}\label{eq: the its vertex of f*(Delta(s,3)) in the case j=2}
f^*((i,m+2) e_{I_s})=
\frac{e_{I_s}}{f^*((i,{m+1}) e_{I_s})} =
\frac{e_{I_s}}{(\sigma(i),m+1) e_{I_s}}
= (\sigma(i),m+2) e_{I_s}\,.
\end{equation}
By \eqref{eq: how it looks like1'} and
\eqref{eq: the its vertex of f*(Delta(s,3)) in the case j=2},
the function $(\sigma(i),m+2) e_{I_s}$ must be the
$i^{\rm th}$ vertex of the simplex
$f^*(\Delta_{s,3})=\sigma' \Delta_{s,j'}$. The latter
shows that $j'=3\ne j$ and $\sigma'(i)=\sigma(i)$ for any $i$
and hence $\sigma'=\sigma$.

Finally, for $j=3$, in the same way as above,
we obtain $j'=2\ne j$ and $\sigma'=\sigma$.
This completes Step 4 and proves Claim.
\vskip7pt

\noindent Continuing the proof of the theorem,
notice that by almost the same argument as in
the Steps 1 and 2 above, one can show that there is
a permutation $\vartheta\in {\mathbf S}(m-1)\subset{\mathbf S}(n)$
and an index $l\in\{1,2,3\}$ such that
\begin{equation}\label{eq: how it looks like2.5}
f^*(\Delta_{s,1}) =\vartheta\Delta_{s,l}\quad\text{for any}\ \
s\ge m+3\,.
\end{equation}
It follows from Claim that either $l=1$ or $l\in\{2,3\}=\{j,j'\}$,
where $j$ and $j'$ are defined by $(ii)$ and $(ii')$ and,
according to Step 4 of the proof above, are distinct elements
of the set $\{2,3\}$.

Thus, we must consider the following three cases:
$a)$ $l=1$, $b)$ $l=j$ and $c)$ $l=j'$.

$a)$ In this case \eqref{eq: how it looks like2.5} takes the form
$f^*(\Delta_{s,1}) = \vartheta\Delta_{s,1}$; therefore,
using \eqref{eq: action of permutation on composition}
for all vertices $\mu$ of the simplex $\Delta_{s,1}$
and the permutation $\rho=\vartheta$, we obtain
\begin{equation}
\label{eq: image of Delta(s,1)}
(\alpha(\vartheta)f)^*(\Delta_{s,1}) =
\vartheta^{-1}[f^*(\Delta_{s,1})]
 = \vartheta^{-1}\vartheta\Delta_{s,1}=\Delta_{s,1}\,;
\end{equation}
in terms of the vertices of the ordered simplices,
this means that for any $i=1,...,m-1$
$$
(\alpha(\vartheta)f)^*((i,m)e_{I_s})=(i,m)e_{I_s}
$$
and
$$
(\alpha(\vartheta)f)^*(e_{I_s})=e_{I_s}
$$
for all $s\ge m+3$. Since $e_{I_{m+3}}$,...,$e_{I_n}$ are all
the vertices of the simplex $\nabla_1^{n-m-3}$, we see that
$$
(\alpha(\vartheta)f)^*
(\nabla_1^{n-m-3})=\nabla_1^{n-m-3}
$$
and
$$
(\alpha(\vartheta)f)^*
((i,m)\nabla_1^{n-m-3})=(i,m)\nabla_1^{n-m-3}
$$
for all $i=1,...,m-1$; this proves the theorem in the case $a)$.

Let us prove now that the cases $b)$ and $c)$ are impossible.
Indeed, if $l=j$ then
\eqref{eq: how it looks like2.5} and $(ii)$ imply that
$f^*(\Delta_{s,1}) = \varphi [f^*(\Delta_{s,2})]$ with
the permutation $\varphi=\vartheta\sigma^{-1}
\in{\mathbf S}(m-1)\subset{\mathbf S}(n)$
that does not depend on $s$. Therefore, using
\eqref{eq: action of permutation on composition additional}
for the vertices $\mu$ of the simplex $\Delta_{s,2}$
and the permutation $\rho=\varphi$, we obtain
$ f^*(\Delta_{s,1}) =  f^*(\alpha(\varphi)\Delta_{s,2})$;
since $f^*$ is an automorphism of the complex
$L_\vartriangle({\mathcal E}^n)$
(see Lemma \ref{Lm: f* isomorphism}),
the latter relation implies
\begin{equation}
\label{eq: equality between permuted simplices Delta1 and Delta 2}
\Delta_{s,1}=\alpha(\varphi)\Delta_{s,2}\,.
\end{equation}
Notice that
$e_{I_s}$ is the very last vertex of the $(m-1)$-simplices
$\Delta_{s,1}$ and $\Delta_{s,2}$; hence for all $s\ge m+3$
we have $e_{I_s}=\alpha(\varphi)e_{I_s}$.
Since $e_{I_{m+3}}$,...,$e_{I_n}$ are all
the vertices of the simplex $\nabla_1^{n-m-3}$, we see that
$\nabla_1^{n-m-3}=\alpha(\varphi)\nabla_1^{n-m-3}$.
Lemma \ref{Lm: stabilizer of nabla1 maxdim} implies
that $\alpha(\varphi)\in {\mathbf S}(m-1)\subset{\mathbf S}(n)$.
As $m>1$, one can see from the definition of
$\Delta_{s,1}$ and $\Delta_{s,2}$ that
$\Delta_{s,1}\ne\psi\Delta_{s,2}$ for any
$\psi\in {\mathbf S}(m-1)$, which contradicts to
\eqref{eq: equality between permuted simplices Delta1 and Delta 2}.
The case $c)$ may be treated similarly. This completes the
proof of the theorem in the case $n>m+3$.

When $n=m+3$, the theorem asserts that 
there is a permutation $\rho\in {\mathbf S}(n)$ such that
$(\rho f)^*(e_{{\mathbf m}(\hat r);r,m+1,m+2,m+3})
=e_{{\mathbf m}(\hat r);r,m+1,m+2,m+3}$ for any
$r\in \{1,...,m\}$. 
The functions $e_{{\mathbf m}(\hat 1);1,m+1,m+2,m+3}$,...,
$e_{{\mathbf m}(\hat m);m,m+1,m+2,m+3}$ are all the vertices
of the ordered simplex $\nabla_2^{m-1}$. Thus,
the statement in question says that
$(\rho f)^*(\nabla_2^{m-1})=\nabla_2^{m-1}$
for an appropriate permutation $\rho\in {\mathbf S}(n)$.
Since $f$ is strictly equivariant,
this is equivalent to the existence of $\phi\in{\mathbf S}(n)$
such that $f^*(\nabla_2^{m-1})=\phi \nabla_2^{m-1}$, which
means precisely that the simplex $f^*(\nabla_2^{m-1})$
is of the second type. The latter property follows from Lemma 
\ref{Lm: f* isomorphism}, which completes the proof.
\end{proof}

\section{The proof of Theorem
\ref{Thm: endomorphisms of En(X)}}
\label{The main theorem}
\noindent Here we prove the main result of this paper.
We start with the following remark which is similar to
Remark 2.14 in \cite{LinSphere}.

\begin{Remark}
\label{Rmk: equality for generic points implies equality everywhere}
For $n\ge m+3$, there exists a non-empty Zariski open subset
$U\subset{\mathcal E}^n(X,gp)$ such that
if $Aq=\sigma q$ for some $q\in U$,
$A\in\mathbf{PSL}(m+1,\mathbb{C})$ and $\sigma\in{\mathbf S}(n)$
then $A=\Id$ and $\sigma =\Id$.

Indeed, Lemma \ref{Lm: existence of the special map in the space}
implies that for any two points $q=(q_1,...,q_n)$ and
$q'=(q'_1,...,q'_n)$ in ${\mathcal E}^n(X,gp)$,
an element $A\in \mathbf{PSL}(m+1,\mathbb{C})$
is uniquely determined by the requirement $Aq_i=q'_i$ for all
$i=1,...,m+2$. Since ${\mathbf S}(n)$ is finite,
it follows that the set $S$ of all points
$q=(q_1,...,q_n)\in{\mathcal E}^n(X,gp)$
such that for some $A\in\mathbf{PSL}(m+1,\mathbb{C})$
and some non-trivial permutation $\sigma\in{\mathbf S}(n)$
(both $A$ and $\sigma$ may depend on $q$) the point
$Aq=(Aq_1,...,Aq_n)$ coincides with the permuted point
$\sigma q=(q_{\sigma^{-1}(1)},...,q_{\sigma^{-1}(n)})$
is a proper Zariski closed subset of ${\mathcal E}^n(X,gp)$.
Its complement $U={\mathcal E}^n(X,gp)\setminus S$
is a desired non-empty Zariski open set.
\hfill $\bigcirc$
\end{Remark}
\vskip7pt

\subsection{Proof of Theorem \ref{Thm: endomorphisms of En(X)}}
\label{ss: Proof of main Theorem}
According to Theorem
\ref{Thm: cross ratio is invariant under automorphism for dim>2},
there exists a permutation $\rho$ such that
\begin{equation}\label{eqv: inv after reordering multidim}
\begin{split}
e_{{\mathbf m}(\hat r);r,m+1,m+2,s}(\rho
f(q))&=e_{{\mathbf m}(\hat r);r,m+1,m+2,s}(q)
\end{split}
\end{equation}
for all $q\in \Cal{E}^n(X,gp)$, $s=m+3,...,n$ and
$r=1,...,m$. Lemma
\ref{Lm: existence of the special map in the space}
implies that there exists a map
$\gamma\colon\Cal{E}^n(\Bbb{CP}^m,gp)\to \mathbf{PSL}(m+1,\mathbb{C})$ such that
$\gamma(q)q\in M_{m,n}$ (see Definition \ref{Def: M subspace}).
Lemma \ref{Lm: invariants of psl(m+1) action} says that
determinant cross ratios are $\mathbf{PSL}(m+1,\mathbb{C})$-invariant;
therefore
\begin{equation}\label{eq: cross ratios are invariant under psl action for q}
e_{{\mathbf m}(\hat r);r,m+1,m+2,s}(\gamma(q)q)
=e_{{\mathbf m}(\hat r);r,m+1,m+2,s}(q)
\end{equation}
and
\begin{equation}\label{eq: cross ratios are invariant under psl action for f(q)}
e_{{\mathbf m}(\hat r);r,m+1,m+2,s}(\gamma(\rho f(q))\rho f(q))
=e_{{\mathbf m}(\hat r);r,m+1,m+2,s}(\rho f(q))
\end{equation}
for all $q\in \Cal{E}^n(X,gp)$,
$s\in \{m+3,\dots,n\}$ and $r\in \{1,\dots,m\}$.
Comparing
\eqref{eqv: inv after reordering multidim},
\eqref{eq: cross ratios are invariant under psl action for q} and
\eqref{eq: cross ratios are invariant under psl action for f(q)}
we obtain that
\begin{equation}\label{eq: the equality for all r and s}
e_{{\mathbf m}(\hat r);r,m+1,m+2,s}(\gamma(q)q)
=e_{{\mathbf m}(\hat r);r,m+1,m+2,s}
(\gamma(\rho f(q))\rho f(q))
\end{equation}
for any $q\in \Cal{E}^n(X,gp)$ and all $s\in \{m+3,\dots,n\}$
and $r\in \{1,\dots,m\}$. Both points
$\gamma(q)q$ and $\gamma(\rho f(q))(\rho f(q))$ are in
$M_{m,n}$, and Lemma \ref{Lm: definition and properties of P}
says that the functions $e_{{\mathbf m}(\hat r);r,m+1,m+2,s}$ with
$s\in \{m+3,\dots,n\}$ and $r\in \{1,\dots,m\}$
separate points of $M_{m,n}$. Consequently,
\eqref{eq: the equality for all r and s} implies that
$\gamma(\rho f(q))\rho f(q) =\gamma(q)q$, or, which is the same,
$\rho f(q) =(\gamma(\rho f(q)))^{-1}\gamma(q)q$.
Set $\tau(q)=(\gamma(\rho f(q)))^{-1}\gamma(q)$;
the map
$$
\tau\colon \Cal{E}^n(X,gp) \ni q
\mapsto \tau(q)\in\mathbf{PSL}(m+1,\mathbb{C})
$$
is  holomorphic and $\tau(q)q=\rho f(q)$, that is,
$f(q)=\sigma\tau(q)q$ for all $q\in{\mathcal E}^n(X,gp)$,
where $\sigma=\rho^{-1}\in{\mathbf S}(n)$.

To complete the proof, we must check that the morphism $\tau$ is
${\mathbf S}(n)$-invariant. Let $\alpha\in\Aut\mathbf{S}(n)$
be the automorphism related to our strictly equivariant map $f$.
For every $\theta\in{\mathbf S}(n)$ and all
$q\in{\mathcal E}^n(X,gp)$ we have
$$
\sigma\tau(\theta q)\theta q
=f(\theta q)=\alpha(\theta)f(q)
=\alpha(\theta)\sigma\tau(q)q\,,
$$
which can be written as
\begin{equation}\label{eq: a}
[(\tau(\theta q))^{-1}\cdot\tau(q)]q
=\sigma^{-1}\alpha(\theta^{-1})\sigma\theta q\,,
\end{equation}
where $(\tau(\theta q))^{-1}\cdot\tau(q)\in
\mathbf{PSL}(m+1,\mathbb{C})$ means
the product of two elements of
the group $\mathbf{PSL}(m+1,\mathbb{C})$. In view of Remark
\ref{Rmk: equality for generic points implies equality everywhere},
this implies that
\begin{equation}
\label{eq: sigma alpha phi equation and tau phi equation}
\sigma^{-1}\alpha(\theta^{-1})\sigma\theta=\Id \ \ \text{and} \  \
\tau(\theta q)=\tau(q)
\end{equation}
for all $\theta\in{\mathbf S}(n)$ and
all $q$ in a non-empty Zariski open subset
$U\subset{\mathcal E}^n(X,gp)$.
Since $\tau$ is continuous, the latter relation holds true for all
$q\in {\mathcal E}^n(X,gp)$; since $\theta\in{\mathbf S}(n)$
was arbitrary, this shows that the morphism
$\tau\colon{\mathcal E}^n(X,gp)\to\mathbf{PSL}(m+1,\mathbb{C})$
is ${\mathbf S}(n)$-invariant. This completes the proof of
Theorem \ref{Thm: endomorphisms of En(X)}.
\hfill $\square$
\vskip7pt

\noindent The following statement is an obvious corollary of
Theorem \ref{Thm: endomorphisms of En(X)}
and Definition \ref{Def: tame map}.

\begin{Corollary}
Let $m>1$, $n\ge m+3$ and $n\ne 2m+2$.

$a)$ Any holomorphic map
$F\colon\Cal{C}^n(\mathbb{CP}^m,gp)\to\Cal{C}^n(\mathbb{CP}^m,gp)$
that can be lifted to a strictly equivariant holomorphic map
$f\colon\Cal{E}^n(\mathbb{CP}^m,gp)\to\Cal{E}^n(\mathbb{CP}^m,gp)$
is tame.

$b)$  Any holomorphic map
$F\colon\Cal{C}^n(\mathbb{C}^m,gp)\to
\Cal{C}^n(\mathbb{C}^m,gp)$ that can be lifted to
a strictly equivariant holomorphic map
$f\colon\Cal{E}^n(\mathbb{C}^m,gp)\to\Cal{E}^n(\mathbb{C}^m,gp)$
is quasitame.
\end{Corollary}

\bibliographystyle{amsplain}

\begin{thebibliography}{1}


%

\bibitem{Barvinok86}
A. I. Barvinok, \emph{Homological type of spaces of configurations of
  structurally stable type in $\mathbb{C}^2$}, Matematicheskie Zametki
  \textbf{Vol 39, No. 1} (1986), 108--112.







\bibitem{Efimov(rus)}
N. V. Efimov, \emph{Visshaya geometriya (advanced geometry)},
Nauka, 1971. English translation: Efimov, N. V. \emph{Higher
geometry.} Translated from the sixth Russian edition by P. C.
Sinha. "Mir", Moscow; distributed by Imported Publications,
Chicago, Ill., 1980.



\bibitem{Feler04a}
Y. Feler,
{\em Configuration spaces of tori},
unpublished, Electronic version arXiv:math.CV/0601283.

\bibitem{GR}
R. C. Gunning, H. Rossi, {\em Analytic functions of several complex
variables}, Prentice-Hill, 1965.







\bibitem{Lin72a} 
V. Ya. Lin,
{\em On representations of the braid group by permutations},
Uspehi Matem. Nauk {\bf 27}, n. 3 (1972), 192.

\bibitem{Lin72b}
\bysame,
{\em Algebraic functions with universal discriminant
manifolds},
Funct. Anal. Appl. \textbf{Vol. 6, No.1} (1972), 73--75.

\bibitem{Lin72c} 
\bysame,
{\em On superpositions of algebraic functions},
Funkt. Analiz i Priloz. {\bf 6}, no. 3 (1972), 77--78
(in Russian). {\em English translation:}
Funct. Anal. Appl. {\bf 6}, n. 3 (1972), 240--241.

\bibitem{Lin79} 
\bysame, {\em Artin braids and the groups and spaces connected
with them}, Itogi Nauki i Tekhniki, Algebra, Topologiya,
Geometriya {\bf 17}, VINITI, Moscow, 1979, pp. 159--227 (in
Russian). {\em English translation:} Journal of Soviet Math.
{\textbf 18} (1982), 736--788.

\bibitem{Lin96b}
\bysame,
{\em Braids, permutations, polynomials}--I,
Preprint MPI 96-118, Max-Planck-Institut f\"ur Mathematik in Bonn,
August 1996, 112pp.

\bibitem{Lin04b}
\bysame,
{\em Braids and permutations},
Electronic version arXiv:math.GR/0404528.

\bibitem{LinSphere}
\bysame, {\em Configuration spaces of
$\mathbb C$ and ${\mathbb{CP}}^1$:
some analytic properties}, Max-Planck-Institut f{\" u}r
Mathematik Preprint Series 2003 (98), Bonn 2003, 80pp.
Revised electronic version arXiv:math.AG/0403120.

\bibitem{Mobius27}
A. F. M{\" o}bius, {\em Der Barycentrishe Calcul},
Verlag von Johann Ambrosius Barth, Leipzig, 1827.

\bibitem{Moulton98}
V. L. Moulton, {\em Vector braids},
   J. Pure Appl. Algebra,
  \textbf{131}, 1998, no.~3, 245--296, \MR{99f:57001}.

\bibitem{Shafarevich94}
I. R. Shafarevich, {\em Basic Algebraic Geometry 1},
Springer-Verlag, 1994.


\bibitem{TerasomaFGM}
T. Terasoma,
{\em Fundamental group of the moduli space of arrangements},
preprint.
Electronic version: http://gauss.ms.u-tokyo.ac.jp/paper/ps/fun.ps.

\bibitem{Weyl46}
H. Weyl, \emph{The {C}lassical {G}roups. {T}heir {I}nvariants and
  {R}epresentations}, Princeton University Press, Princeton, N.J., 1946.


\bibitem{Zinde74}
V.~M. Zinde, {\em Commutants of Artin groups}, Uspehi Mat. Nauk \textbf{30} (1975), no.~5(185), 207--208. \MR{55 \#5750}.

\bibitem{Zinde77a}
\bysame, {\em Analytic properties of the spaces of regular orbits of {C}oxeter groups of the series {$B$} and {$D$}}, Funkcional. Anal. i Prilo\v zen. \textbf{11} (1977), no.~1, 69--70. \MR{58 \#1243}.

\bibitem{Zinde77b}
\bysame, {\em Holomorphic mappings of the spaces of regular orbits of
Coxeter groups of series {$B$} and {$D$}}, Sibirsk. Mat. \v Z. \textbf{18} (1977), no.~5, 1015--1026, 1205. \MR{57 \#706}.

\bibitem{Zinde77c}
\bysame, {\em Some homomorphisms of the Artin groups of the series
{$B\sb{n}$} and {$D\sb{n}$} into groups of the same series and into symmetric groups}, Uspehi Mat. Nauk \textbf{32} (1977), no.~1(193), 189--190. \MR{56 \#3129}.

\bibitem{Zinde78}
\bysame, {\em Studies of homomorphisms of Artin groups}, C. R. Math. Rep. Acad. Sci. Canada \textbf{1} (1978/79), no.~4, 199--200. \MR{80j:20039}.

\end{thebibliography}

\end{document}